\documentstyle[amssymb,amsfonts]{amsart}

\font \roman = cmr10 at 10 true pt

\def\R{{\hbox{\bf R}}}

\def\Q{{\hbox{\bf Q}}}
\def\B{{\hbox{\bf B}}}

\def\f{{\overline{f}}}
\def\x{{\hbox{\roman x}}}
\def\y{{\hbox{\roman y}}}
\def\z{{\hbox{\roman z}}}

\def\E{{\hbox{\bf E}}}

\def\T{{\hbox{\bf T}}}

\def\td{{T_\delta}}
\def\Td{{\T_\delta}}
\def\dt{{N\delta}}
\def\tdt{{T_\dt}}
\def\Tdt{{\T_\dt}}
\def\sd{{\rho}}
\def\sdt{{N\rho}}
\def\se{{\sqrt{\epsilon}}}

\def\tsd{{T_\sd}}
\def\tsdo{{T_\sd^0}}
\def\tsda{{T_\sd^1}}
\def\tsdb{{T_\sd^2}}
\def\tsdc{{T_\sd^3}}
\def\tsdd{{T_\sd^4}}
\def\tsdk{{T_\sd^k}}
\def\tsdka{{T_\sd^{k+1}}}
\def\tsdn{{T_\sd^{n-1}}}
\def\tsdnn{{T_\sd^n}}
\def\tsdi{{T_\sd^i}}
\def\Tsd{{\T_\sd}}

\def\mult{{\mu}}
\def\allt#1{%
\smash{
\vtop{%
     \ialign{%
        ##\crcr
        $\hfil\displaystyle{\tilde \forall}\hfil$\crcr%
        \noalign{\kern1.5pt\nointerlineskip}
        $\hfil\!\!#1\hfil$\crcr\noalign{\kern1.5pt}
        }
       }
      } \hbox{$\vphantom{#1}$}
     }

\def\be#1{\begin{equation}\label{#1}}
\def\bas{\begin{align*}}
\def\eas{\end{align*}}
\def\bi{\begin{itemize}}
\def\ei{\end{itemize}}

\def\dir{{\hbox{\roman dir}}}
\def\dimsup{\overline{{\hbox{\roman dim}}}}
\def\dist{{\hbox{\roman dist}}}

\def\eps{\varepsilon}

\newenvironment{proof}{\noindent {\bf Proof} }{\endprf\par}
\def \endprf{\hfill  {\vrule height6pt width6pt depth0pt}\medskip}
\def\emph#1{{\it #1}}
\def\textbf#1{{\bf #1}}

\def\Ecal{{\cal E}}

\def\N{{\cal N}}

\theoremstyle{plain}
  \newtheorem{theorem}[subsection]{Theorem}

  \newtheorem{proposition}[subsection]{Proposition}
  
  \newtheorem{lemma}[subsection]{Lemma}

\theoremstyle{remark}

\theoremstyle{definition}
  \newtheorem{definition}[subsection]{Definition}

\include{psfig}

\begin{document}

\title[Besicovitch sets in medium dimension]{An improved bound for the Minkowski dimension of Besicovitch sets in medium dimension}

\author{Izabella {\L}aba}
\address{Department of Mathematics, Princeton University, Princeton NJ
08544}
\email{laba@@math.princeton.edu}

\author{Terence Tao}
\address{Department of Mathematics, UCLA, Los Angeles CA 90095-1555}
\email{tao@@math.ucla.edu}


\begin{abstract}  We use geometrical combinatorics arguments, including the
``hairbrush'' argument of Wolff \cite{W1}, the x-ray estimates in \cite{W2},
\cite{LT1}, and the sticky/plany/grainy analysis of \cite{KLT}, to show 
that Besicovitch sets in $\R^n$ have Minkowski dimension at least 
$\frac{n+2}{2} + \eps_n$ for all $n \geq 4$, where $\eps_n > 0$ is an 
absolute constant depending only on $n$.  This complements the results of 
\cite{KLT}, which established the same result for $n=3$, and of \cite{B3},
\cite{KT}, which used arithmetic combinatorics techniques to establish the
result for $n \geq 9$.  Unlike the arguments in \cite{KLT}, \cite{B3}, 
\cite{KT}, our arguments will be purely geometric and do not require 
arithmetic combinatorics.
\end{abstract}

\maketitle

\section{Introduction}

Let $n \geq 2$ be an integer.  We recall the following definitions:

\begin{definition}  A \emph{Besicovitch set} (or ``Kakeya set'')
$E \subset \R^n$ is a set
which contains a unit line segment in every direction.
\end{definition}

\begin{definition}\label{mink-def}  For any set $E \subset \R^n$, the 
(upper) \emph{Minkowski dimension} $\dimsup(E)$ is defined as
$$ \dimsup(E) = n - \liminf_{\delta \to 0} \log_\delta |\N_\delta(E)|.$$
Here and in the sequel, $\N_\delta(E)$ denotes the $\delta$-neighbourhood of $E$.
\end{definition}
 
Informally, the Kakeya conjecture (see e.g. \cite{B1}) states that all 
Besicovitch sets in $\R^n$ have full dimension; this conjecture has been 
verified for $n = 2$ but is open otherwise.  For the purposes of this paper
we shall restrict ourselves to the upper Minkowski dimension; the 
corresponding problems for Hausdorff dimension or lower Minkowski dimension
are more difficult, and do not seem to be easily attacked by the techniques
in this paper (see the discussion in \cite{KLT}, Section 1).

We briefly summarize some recent progress on this problem.  For a more 
thorough treatment of these problems and their applications see \cite{B1}, 
\cite{wolff:survey}, \cite{wolff:icm}.

In $\R^n$, Wolff \cite{W1} used geometric combinatorics techniques, including
the construction of ``hair-brushes'', to show the estimate
\be{wolff-dim}
 \dimsup(E) \geq \frac{1}{2}n + 1.
\end{equation}
More recently, a very different approach of Bourgain \cite{B3} based on the 
arithmetic combinatorics of Gowers \cite{gowers}, and then developed further
by Katz and Tao \cite{KT} has shown 
\be{borg-est}
 \dimsup(E) \geq \frac{4}{7}n + \frac{3}{7}.
\end{equation}
This improves on \eqref{wolff-dim} when $n \geq 9$.  

By combining the arithmetic combinatorics techniques in \cite{B3} with 
geometric arguments, in particular Wolff's x-ray estimate \cite{W2} and the
observations that a minimal-dimension Besicovitch set must be ``sticky'', 
``plany'', and ``grainy'', Katz, \L aba, and Tao \cite{KLT} managed to 
obtain a small improvement to \eqref{wolff-dim} in the $n=3$ case, namely
$$ 
\dimsup(E) \geq \frac{5}{2} + \eps_3
$$
for some absolute constant $\eps_3 > 0$ ($\eps_3 = 10^{-10}$ will suffice).

The purpose of this paper is to show a similar estimate in higher dimensions:

\begin{theorem}\label{main-thm}
For all $n \geq 4$ and Besicovitch sets $E \subset \R^n$ we have
\be{main-bound}
\dimsup(E) \geq \frac{n+2}{2} + \eps_n
\end{equation}
where $\eps_n > 0$ is an absolute constant depending only on $n$.
\end{theorem}

The bound \eqref{main-bound} is thus already known for $n=3$ and
$n \geq 9$, and is new for $4 \leq n \leq 8$.  We do not attempt to
obtain an optimal value for $\eps_n$, but $\eps_n = (2n)^{-10}$ would certainly suffice.

The arguments of this paper are closely based on those in \cite{KLT}, in 
that they require an x-ray estimate (we shall use the one in \cite{LT1}), 
and the observations of stickiness, planiness, and graininess.  In fact, we 
shall borrow many definitions and lemmas from \cite{KLT} without any 
modifications (other than changing $3$ to $n$ in the obvious places).  
However, the arguments are somewhat simpler than in the $n=3$ case in that 
one does not need to involve arithmetic combinatorial techniques as in 
\cite{B3}, \cite{KT}.  
In fact, the proof is even simpler in the $n > 4$ case, 
mostly because any hypothetical counterexample to \eqref{main-bound} for
$n>4$ would have codimension strictly greater than 1.
Unfortunately, our arguments do not lead to any 
substantial simplifications for the $n=3$ argument in \cite{KLT}, in which 
the codimension is $1/2$.

The remainder of this paper is devoted to a proof of Theorem \ref{main-thm},
after some notational preliminaries in Section \ref{notation}, and is 
organized as follows.
We fix $n \geq 4$, and assume for contradiction that there is a
Besicovitch set $E$ in $\R^n$ with upper Minkowski dimension extremely close to $\frac{n+2}{2}$.

For any scale $0 < \sigma \ll 1$, the $\sigma$-neighbourhood $E_\sigma$ is 
essentially the union of about $\sigma^{1-n}$ tubes with dimensions $\sigma 
\times 1$ and oriented in a $\sigma$-separated set of directions, filling 
out a set of size about $\sigma^{(n-2)/n}$.  We now invoke the x-ray 
estimate in \cite{LT1} (which is a higher-dimensional analogue of Wolff's 
x-ray estimate in \cite{W2}) and the arguments of \cite{KLT} (see also the 
discussion in \cite{W2}) to conclude a certain ``stickiness'' property of 
these tubes in Section \ref{sticky-sec}.  Essentially, this states that if 
$0 < \delta \ll \sigma \ll 1$ and two $\delta$-tubes in $E_\delta$ have 
directions separated by $\lesssim \sigma$, then with high probability these 
two $\delta$-tubes are contained inside a single $\sigma$-tube in $E_\sigma$.
From this stickiness property, and an application of Wolff's Kakeya estimate
\cite{W1} (for instance) at several scales, we can deduce various self-similarity 
properties of $E$ in Section \ref{selfsim-sec}, which informally state that 
various small portions of $E$ have roughly the same size and shape as $E$ 
itself when rescaled appropriately.  This part of the argument is identical 
to that in \cite{KLT}, but generalized to arbitrary dimension.

As in \cite{KLT}, we now analyze $E$ simultaneously at two scales $\delta$ and $\rho$, where $\rho = \sqrt{\delta}$ and $\delta$ is small.  This particular choice of scales has been exploited for many related problems, notably the restriction problem; see e.g. \cite{B1}, \cite{B2}, \cite{TVV}.

The next step in Section \ref{triple-sec} is to deduce a certain ``planiness'' 
property of the $\rho$-tubes in $E_\rho$; roughly speaking, this states that
the $\rho$-tubes that pass through a given point are not spread arbitrarily 
in space, but must be somewhat degenerate.  This follows the philosophy of 
\cite{KLT}, although our notion of degeneracy is slightly different in higher 
dimensions than in the $n=3$ case.
In the high-dimensional case $n > 4$ one can actually show that the $\rho$-tubes 
through a point $x$ must lie in a small neighbourhood of a space of codimension 
at least 2; this follows from the ``planiness-graininess'' relationship in 
\cite{KLT} and the observation that $E$ has co-dimension strictly  greater 
than 1 when $n > 4$.  This leads to a fairly simple way to improve 
\eqref{wolff-dim}, which we pursue in Section \ref{planar-sec}. 
Basically, we continue Wolff's ``hairbrush" argument \cite{W1} and prove
that any two hairbrushes with intersecting stems are either largely disjoint
or concentrated in a small neighbourhood of a space of codimension at least 
1; both of these cases are then easily handled.  

The only case remaining is when $n=4$ and the $\rho$-tubes that pass through 
a given point lie in a small neighbourhood of a genuinely three-dimensional 
object such as a hyperplane.  In this case we again use the ``planiness-graininess'' 
relationship of \cite{KLT}, and deduce that $E_\delta$ has a very specific 
structure locally.  In fact, when localized to balls of radius slightly 
larger than $\delta$, the set $E_\delta$ must look like the $\delta$-neighbourhood 
of a hyperplane.  One can now obtain a gain to \eqref{wolff-dim} arising from 
the geometric fact that there are only a restricted number of possible directions 
of $\delta$-tubes which can pass through four distinct
$\delta$-neighbourhoods of hyperplanes at four separated places (cf. the use 
of the ``three-line lemma'' in \cite{schlag:kakeya}).  We perform this in Section \ref{4D-sec} and 11.

The second author is supported by grants from the Sloan and Packard foundations.

\section{Notation and preliminaries I.}\label{notation}

We shall stay as close to the notation of \cite{KLT} as possible, though of course we are no longer working in $\R^3$.

Throughout this paper $n \geq 4$ will be fixed, with all constants implicitly 
depending on $n$.  We shall fix $d = \frac{n+2}{2}$; this is the lower bound 
on the dimension of Kakeya sets given by \eqref{wolff-dim}.

We shall always be working in $\R^n$.  We use italic letters $x,y,z$
to denote points in $\R^3$, and $\x_i, \y_i, \z_i$ to denote their co-ordinates, 
thus $x = (\x_1, \ldots, \x_n)$.

Unless otherwise specified, all integrals will be over  $\R^n$ with
Lebesgue measure.

In this paper $\delta$ refers to a number such that $0 < \delta \ll 1$,
and $\eps$ refers to a fixed number such that $0 < \eps \ll 1$.  In addition 
to the scale $\delta$, we shall need the intermediate scales 
$$ \delta \ll N\delta \ll \sd \ll N\sd \ll 1$$
defined by
$$ \sd := \sqrt{\delta}; \quad N = \delta^{-K^{-1}}$$
where $K$ is an absolute constant depending only on $n$ ($K = 10n$ will do).

We use $C$, $c$ to denote generic positive constants, varying
from line to line (unless subscripted), which are 
independent of $\eps$, $\delta$, $K$, but which may depend on $d$, $n$.
$C$ will denote the large constants
and $c$ will denote the small constants.  

We will use $X \lesssim Y$, $Y \gtrsim X$, or $X = O(Y)$ to denote
the inequality $|X| \leq A Y$, where $A$ is a positive quantity
which may depend on $\eps$.  We use
$X \gg Y$ to denote the statement $X \geq AY$ for a large constant $A$.
We use $X \sim Y$ to denote the statement that $X \lesssim Y$ and
$Y \lesssim X$.  

We will use $X \lessapprox Y$, $Y \gtrapprox X$, or ``$Y$ majorizes $X$''
to denote the inequality
$$ |X| \leq A \delta^{-C\eps} Y,$$
where $A$ is a positive quantity which may depend on $\eps$, and $C$
is a quantity which does not depend on $\eps$.
We use $X \approx Y$ to denote the statement that $X \lessapprox Y$
and $Y \lessapprox X$.  In particular we have $\eps \approx 1$.

If $E$ is a subset of $\R^n$, we use $|E|$ to denote its Lebesgue measure;
if $I$ is a finite set, we use $\# I$ to denote its cardinality.

As in \cite{KLT}, \cite{LT1}, it will be convenient to define $\sigma$-tube in an ``affine'' manner.  Namely, for any $\delta \leq \sigma \leq 1$
we define a \emph{$\sigma$-tube} $T_\sigma$ to be a $\sigma$-neighbourhood 
of a line segment 
whose endpoints $x$ and $y$ are on the planes $\{ \x_n = 0\}$ and
$\{ \y_n = 1\}$ respectively, and whose orientation is within
$\frac{1}{10}$ of the vertical.  We call $y-x$ the \emph{direction} of 
$T_\sigma$ and denote it by $\dir(T_\sigma)$.  We call $\sigma$ the 
\emph{thickness} of $T_\sigma$.  Whenever possible, we shall try to 
subscript a tube by its thickness.  Note that
\be{ts-volume}
 | T_\sigma | \sim \sigma^{n-1}
\end{equation}
for any $\sigma$-tube $T_\sigma$.  

If $T$ is a tube, we define $CT$ to be the dilate of $T$
about its axis by a factor $C$.   We say that two tubes $T$ and $T'$
are \emph{equivalent} if $T \subset CT'$ and $T' \subset CT$.
If $\T$ is a set of tubes, we say that $\T$ \emph{consists of essentially
distinct tubes} if for any $T \in \T$ there are at most $O(1)$ tubes $T'$
which are equivalent to $T$.

We use the term $r$-ball to denote a ball of radius $r$, and
use $B(x,r)$ to denote the $r$-ball centered at $x$.  

If $1 < p < \infty$ is an exponent, we define the dual exponent
by $p' = p/(p-1)$.

\section{X-ray estimates}\label{xray}

In this section we summarize the x-ray estimate from \cite{LT1} which
we shall need, especially in the proofs of Propositions \ref{sticky-reduction} 
and \ref{p3-refine}.  In the following $\sigma, \theta$ are quantities
such that $\delta \leq \sigma \leq \theta \leq 1$.

\begin{definition}  If $\T_\sigma$ is a collection of $\sigma$-tubes, we 
define the \emph{directional multiplicity} $m = m(\T_\sigma)$ to be the largest
number of tubes in $\T_\sigma$ whose directions all lie in a 
cap of radius $\sigma$.  If $m \approx 1$, we say that
$\T_\sigma$ is \emph{direction-separated}.
\end{definition}

\begin{lemma}\label{narrow} Let $\delta \leq \sigma \leq \theta \ll 1$, and
let $\T_\sigma$ be a collection of essentially distinct
$\sigma$-tubes with directional multiplicity at most $m$, and whose set of 
directions all lie in a $\sigma$-separated set $\E$.  Then we have
\be{max-xray}
 \| \sum_{T_\sigma \in \T_\sigma} \chi_{T_\sigma} \|_{d'}
\lessapprox \sigma^{\frac{d-n}{d}} (\sigma^{n-1} \# \E)^{\frac{n-2}{n-1} + 
\frac{1}{d(n-1)}} 
m^{1-\beta},
\end{equation}
where $\beta > 0$ is an absolute constant depending only on $n$.  In 
particular, if $\E$ is contained in a cap of width $\theta$, then
$$
\| \sum_{T_\sigma \in \T_\sigma} \chi_{T_\sigma} \|_{d'}
\lessapprox \sigma^{\frac{d-n}{d}} \theta^{n-2 + \frac{1}{d}}
m^{1-\beta}.
$$
\end{lemma}
 
In the notation of \cite{KLT}, we are stating that we have an x-ray estimate 
at dimension $d = \frac{n+2}{2}$ in $\R^n$.  

\begin{proof}  
By a direct application of \cite{LT1}, Theorem 1.2 we have
$$
\| \sum_{T_\sigma \in \T_\sigma} \chi_{T_\sigma} \|_{d'}
\lessapprox \sigma^{\frac{d-n}{d}} m^{1/q - 1/r} (\sigma^{n-1} |\T_\sigma|)^{1/q'},
$$
where $q = \frac{(n-1)(n+2)}{n}$, $r = 2(n+2)$.  From the assumptions on 
$\T_\sigma$ we see that
$$ |\T_\sigma| \lessapprox m \# \E,$$
and the claim follows by some algebra (with $\beta = \frac{1}{2(n+2)}$).
\end{proof}

The result in \cite{LT1} is a higher-dimensional version of the x-ray 
estimate in \cite{W2}.  We remark that if one only had the Kakeya estimates 
from \cite{W1} available then one could only show \eqref{max-xray} with 
$\beta = 0$.

\section{The sticky reduction}\label{sticky-sec}

We now begin the proof of Theorem \ref{main-thm}.  We assume for contradiction 
that there exist Besicovitch sets of upper Minkowski dimension at most 
$d+\eps$.  We shall eventually show that this leads to a contradiction if 
$\eps$ was sufficiently small.

As in \cite{KLT}, we shall use the hypothesis of a Besicovitch set of 
near-minimal Minkowski dimension to obtain a ``sticky'' collection of tubes,
which we now pause to define. 

\begin{definition}  Let $\Td$ be a collection of $\delta$-tubes.
We say that $\Td$ is \emph{sticky} if 
$\Td$ is direction-separated and there exists
a collection $\Tsd$ of direction-separated $\sd$-tubes 
and a partition of $\Td$ into disjoint sets $\Td[\tsd]$ for $\tsd \in \Tsd$ such that
\be{sticky-max}
\td \subset \tsd \hbox{ for all } \tsd \in \Tsd \hbox{ and } \td \in \Td[\tsd],
\end{equation}
and we have the cardinality estimates
\begin{align}
\# \Td &\approx (\frac{1}{\delta})^{n-1}\label{t-card}\\
\# \Tsd &\approx (\frac{1}{\sd})^{n-1}\label{tsd-card}\\
\# \Td[\tsd] &\approx (\frac{\sd}{\delta})^{n-1} \hbox{ for all } \tsd \in \Tsd.
\label{filled}
\end{align}
\end{definition}

The following Proposition, which is a consequence of Lemma \ref{narrow}, 
was essentially proven in \cite{KLT}; the extension from three dimensions 
to general dimension is trivial.  (The connection between x-ray estimates 
and the stickiness of Besicovitch sets of near-minimal Minkowski dimension 
was first noted in \cite{W2}).   We shall use a similar argument in the 
proof of Lemma \ref{sticky-dt}.

\begin{proposition}\label{sticky-reduction}\cite{KLT}  Suppose there exists 
a Besicovitch set $E$ with $\dimsup(E) < d + \eps$.
Then for any sufficiently small $\delta$, there exists a 
sticky collection $\Td$ of tubes at scale $\delta$, with the associated
collection $\Tsd$, such that
\be{emush-sd}
|\bigcup_{\tsd \in \Tsd} \tsd| \lessapprox \sd^{n-d}.
\end{equation}
More generally, we have
\be{emush-sigma}
| \N_\sigma(\bigcup_{\td \in \Td} \td)| \lessapprox \sigma^{n-d}
\end{equation}
for all $\delta \leq \sigma \leq 1$.
\end{proposition}

One can obtain stickiness for scales other than $\sd$, but we shall not do 
so here.  (Later on we shall implicitly repeat a version of the sticky 
reduction at scale $\dt$; see the proof of Lemma \ref{sticky-dt}).

\section{More notation}

In the rest of the paper $\Td$ will be a sticky collection of tubes
satisfying \eqref{emush-sd}, \eqref{emush-sigma}.

For future reference we shall set out some notation and estimates which
we shall use frequently.

\begin{definition}  For any $x \in \R^3$ and $\tsd \in \Tsd$, we define the sets
$\Td(x)$, $\Tsd(x)$, and $\T[\tsd](x)$ by
\begin{align*}
\Td(x) &= \{ \td \in \Td: x \in \td\}\\
\Tsd(x) &= \{ \tsd \in \Tsd: x \in \Tsd\}\\
\Td[\tsd](x) &= \{ \td \in \Td[\tsd]: x \in \td\}.
\end{align*}
\end{definition}

\begin{definition}  We define the sets $E_\delta$, $E_\sd$, and
$E_\delta[\tsd]$ for all $\tsd \in \Tsd$ by
\begin{align*}
E_\delta &:= \bigcup_{\td \in \Td} \td\\
E_\sd &:= \bigcup_{\tsd \in \Tsd} \tsd\\
E_\delta[\tsd] &:= \bigcup_{\td \in \Td[\tsd]} \td.
\end{align*}
We similarly 
define the multiplicity functions $\mult_\delta$, $\mult_\sd$, and
$\mult_\delta[\tsd]$ by
\begin{align*}
\mult_\delta(x) &:= \sum_{\td \in \Td} \chi_\td(x) = \# \Td(x) \\
\mult_\sd(x) &:= \sum_{\tsd \in \Tsd} \chi_\tsd(x) = \# \Tsd(x) \\
\mult_\delta[\tsd](x) &:= \sum_{\td \in \Td[\tsd]} \chi_\tsd(x)
= \# \Td[\tsd](x).
\end{align*}
\end{definition}

We borrow the following notation from \cite{KLT}.

\begin{definition}\cite{KLT}  Let $P(x)$ and $Q(x)$ be 
logical statements with free parameters $x = (x_1, \ldots, x_n)$, where 
each of the variables $x_i$ range either over a subset of Euclidean space, 
or over a discrete set.  We use
\be{pq}
\allt{x} Q(x): P(x)
\end{equation}
to denote the statement that
\be{rarity} 
|\{ x: Q(x) \hbox{ holds, but } 
P(x) \hbox{ fails} \}|
\lessapprox \delta^{c\se} 
|\{ x: Q(x) \hbox{ holds}\}|
\end{equation}
for some absolute constant $c > 0$, where the sets are measured
with respect to the measure $dx = \prod_{i=1}^n x_i$, and
$dx_i$ is Lebesgue measure if the $x_i$ range over a subset of
Euclidean space, or counting measure if they range over a discrete set.
\end{definition}

In practice our variables $x_i$ will either be points in $\R^3$
(and thus endowed with Lebesgue measure), or
tubes in $\Td$ or $\Tsd$ (and thus endowed with counting measure).  
Thus, for instance,
$$ \allt{\td,x} \td \in \Td, x \in \td: P(x,\td)$$
denotes the statement that
$$ \sum_{\td \in \Td} | \{ x \in \td: P(x,\td) \hbox{ fails} \} | \lessapprox 
\delta^{c\se} \sum_{\td \in \Td} |\td|.
$$
The
right-hand side of \eqref{rarity} will always be automatically finite
in our applications.  Note that \eqref{pq}
vacuously holds if $P(x)$ is never satisfied.

The statement \eqref{pq} should be read as ``for most $x$
satisfying $P(x)$, $Q(x)$ holds'', where ``most'' means that
the event occurs with probability very close to 1.

We recall the following properties of $\tilde \forall$ from \cite{KLT}.

\begin{lemma}\label{allt-mix}\cite{KLT}
Suppose that $Q_1(x)$, $Q_2(x,y)$, and $P(x,y)$ are properties
depending on some free parameters $x = (x_1, \ldots, x_n)$,
$y = (y_1, \ldots, y_m)$ which obey 
\be{uniform-q2}
| \{ y: Q_2(x,y) \} | \approx M
\hbox{ whenever } Q_1(x) \hbox{ holds}
\end{equation}
for some quantity $M$ independent of $x$.  Then, the statements
\be{q1q2-hyp}
\allt{x,y} Q_1(x), Q_2(x,y): P(x,y)
\end{equation}
and
\be{q1q2-conc}
\allt{x} Q_1(x) : [\allt{y} Q_2(x,y) \hbox{ holds } : P(x,y)]
\end{equation}
are equivalent (up to changes of constants).
\end{lemma}

Here and in the rest of the paper, the expression ``$Q(x),P(x)$''
is an abbreviation for ``$Q(x)$ and $P(x)$ both hold''.

\begin{lemma}\label{delta-sd}\cite{KLT}  Let $\Td$ be a sticky collection of
tubes, and let $P(y,\tsd,\td)$ be a property.  
Then the statements
\be{dsd-hyp}
\allt{\tsd,\td,y} \tsd \in \Tsd, \td \in \Td[\tsd], y \in \td: P(y,\tsd,\td)
\end{equation}
\be{dsd-dum}
\allt{\tsd,\td,y,x} \tsd \in \Tsd, \td \in \Td[\tsd], y \in \td, x \in \tsd \cap
B(y,C\sd):  P(y,\tsd,\td)
\end{equation}
and
\be{dsd-conc}
\allt{\tsd,x} \tsd \in \Tsd, x \in \tsd: 
[ \allt{\td,y} \td \in \Td[\tsd], y \in \td \cap B(x,C\sd): P(y,\tsd,\td) ]
\end{equation}
are equivalent (up to changes of constants).
\end{lemma}

\section{Uniformity and self-similarity}\label{selfsim-sec}

We continue the strategy of \cite{KLT}, and use the stickiness of $\Td$ to 
imply certain self-similarity properties of the set $E_\delta$; roughly 
speaking, we wish to prove a rigorous version of \cite{KLT}, Heuristic 6.1 
with the obvious modifications to $n$ dimensions.  These properties will 
have many uses, but are especially important for deriving planiness and 
graininess properties, as we shall see.

We shall need the following rather technical definitions from \cite{KLT}.

\begin{definition}\label{p3a-def}\cite{KLT}
If $x_0 \in \R^n$ and $\tsd \in \Tsd[x_0]$, we say that
$P_1(x_0,\tsd)$ holds if the three statements 
\be{edim-sd-lower} 
|E_\delta[\tsd] \cap B(x_0,C\sd)| \gtrapprox \delta^{\se}
\delta^{n-d} \sd^d
\end{equation}
\be{edim-uniform}
\allt{\td,x} \td \in \Td[\tsd], x \in \td \cap B(x_0,C\sd): \eqref{edim} 
\hbox{ holds for all } \delta \leq \sigma \ll 1
\end{equation}
\be{enarrow-uniform}
\allt{\td,x} \td \in \Td[\tsd], x \in \td \cap B(x_0,C\sd): \eqref{narrow-local} 
\hbox{ holds}. 
\end{equation}
hold, where \eqref{narrow-local} is the estimate
\be{narrow-local}
\mult_\delta[\tsd](x) \lessapprox \delta^{-\se} \delta^{-(n-d)/2}
\end{equation}
and \eqref{edim} is the estimate
\be{edim}
|E_\delta \cap B(x,C\sigma)| \lessapprox \delta^{-\se} 
\delta^{n-d} \sigma^d.
\end{equation}
\end{definition}

The property \eqref{edim-uniform} states that $E_\delta$ looks locally like 
the $\delta$-neighbourhood of a set of dimension $\leq d$.  The properties 
\eqref{edim-sd-lower} complements these upper bounds on $E_\delta$ by a 
similar lower bound on the individual sets $E_\delta[\tsd]$, while 
\eqref{enarrow-uniform} limits the multiplicity of the tubes $\td$ in 
$\T[\tsd]$.

\begin{definition}\label{p3-def}  
Let $x_0$ be a point in $\R^n$.
We say that $P_2(x_0)$ holds if one has
\begin{align}
\delta^{C\se} \delta^{-(n-d)/2} 
\lessapprox \# \Tsd(x_0) &\lessapprox \delta^{-C\se}
\delta^{-(n-d)/2}
\label{many-fat}\\
\allt{\tsd} \tsd \in \Tsd(x_0)&: P_1(x_0,\tsd) \hbox{ holds}
\label{few-excep}\\
|E_\delta \cap B(x_0,C\sd)| &\lessapprox \delta^{-C\se}
\delta^{n-d} \sd^d\hbox{, and}
\label{edim-sd}\\ 
 \# \{ \tsd \in \Tsd(x_0): \dir(\tsd) \in B(\omega,\theta) \}
&\lessapprox \delta^{-C\se} \theta^c \delta^{-(n-d)/2}
\label{no-cluster}
\end{align}
for all directions $\omega$ and all
$\delta \leq \theta \ll 1$.
\end{definition}

The property \eqref{few-excep} asserts that $P_2$ contains $P_1$ in a 
certain sense. The property \eqref{no-cluster} basically states that
the tubes in $\Tsd(x_0)$ are not clustered in a narrow angular band.  
\eqref{edim-sd} is essentially a re-iteration of \eqref{edim}, while 
\eqref{many-fat} asserts that $x_0$ is contained in the expected number of 
$\sd$-tubes in $\Tsd$.

The following Proposition was essentially proven in \cite{KLT}, with the 
obvious modifications for $\R^n$ (basically, replace any occurrence of the 
number $3$ by $n$ in the proofs of \cite{KLT}, Propositions 6.2, 6.4, 6.6):

\begin{proposition}\label{p3-refine}\cite{KLT}
If the constants in the above definitions are chosen appropriately,
then we have
\be{p1-ok} \allt{\tsd,x} \tsd \in \Tsd, x \in \tsd : P_1(x,\tsd) \hbox{ holds}.
\end{equation}
and
\be{p2-ok} \allt{\tsd,x} \tsd \in \Tsd, x \in \tsd : P_2(x) \hbox{ holds}.
\end{equation}
\end{proposition}

\section{$n-1$-fold intersections}\label{triple-sec}

In this section we fix $x_0$ to be a point in $\R^n$ such that 
$P_2(x_0)$ holds.

Let $\tsd$ be a tube in $\Tsd(x_0)$
such that $P_1(x_0,\tsd)$ holds.  By Proposition
\ref{p3-refine}, this situation occurs almost always.  

Let $A(x_0,\tsd)$ denote the set
\be{ax0tsd-def}
 A(x_0,\tsd) = E_\delta[\tsd] \cap B(x_0,C\sd) = \bigcup_{\td \in \Td[\tsd]}
 \td \cap B(x_0,C\sd).
\end{equation}
From \eqref{edim-sd-lower} we have
$$ |A(x_0,\tsd)| \gtrapprox \delta^{C\se} \delta^{n-d} \sd^d.$$
On the other hand, from \eqref{edim-sd} we have
$$ |\bigcup_{\tsd \in \Tsd(x_0)} A(x_0,\tsd)| 
\lessapprox \delta^{-C\se} \delta^{n-d} \sd^d.$$
Thus we expect a lot of overlap between the $A(x_0,\tsd)$.  In particular,
we expect the size of the $n-1$-fold intersection
\be{intersection}
\bigcap_{i=1}^{n-1} A(x_0,\tsdi) 
\end{equation}
to be quite large for many $n-1$-tuples of tubes $\tsdi$ in $\Tsd(x_0)$.

However, it turns out that we can get a non-trivial bound on the
size of \eqref{intersection} if the tubes $\tsdi)$ are not ``coplanar'', or 
if the sets $A(x_0,\tsdi)$ are not ``grainy'', in the sense that they do 
not resemble the unions of $\delta \times \dt \times \ldots \times \dt$ 
boxes.  (This observation is slightly different from the corresponding 
observation in \cite{KLT}, which dealt with the intersections of $n$ sets 
rather than $n-1$).

More precisely, we have

\begin{definition}  We define a
\emph{square} to be any
rectangular box $Q$ of dimensions $\delta \times \dt \times \ldots \dt$.  
We call the sides of length $\dt$ the \emph{long sides} of $Q$, and we call 
the hyperplane generated by the long sides the \emph{hyperplane} of $Q$.  
We say that $Q$ is \emph{parallel} to a direction $v$ if $v$ is parallel to 
the hyperplane of $Q$.
\end{definition}

\begin{definition}  We say that an $n-1$-tuple $(\tsda, \ldots, \tsdn)$ of 
tubes in $\Tsd$ with a common point $x_0$ is \emph{coplanar} if there 
exists an affine subspace $\Gamma$ of $\R^n$ dimension $n-2$ such that
$$ \tsdi \subset \N_\sdt(\Gamma) $$
for all $1 \leq i \leq n-1$. 
\end{definition}

\begin{lemma}\label{triple}
Let $x_0$ be a point in $\R^n$ such that
$P_2(x_0)$ holds, and let $F$ be a subset of $\R^n$.  Let $\tsda, \ldots, 
\tsdn$ be tubes in $\Tsd(x_0)$ which are not coplanar.  Then
\be{grainy-intersect}
|F \cap \bigcap_{i=1}^{n-1} A(x_0,\tsdi)|
\lessapprox \delta^{-C\se} \delta^{n-d} \sd^d 
(\sup_Q \frac{|F \cap Q|}{|Q|})^{1/n},
\end{equation}
where $Q$ ranges over all
squares parallel to $\dir(\tsda)$.
\end{lemma}

Far stronger versions of this lemma are possible (e.g. one can obtain 
analogues to Lemma 7.3 in \cite{KLT}); however, this form of the Lemma is 
adequate for our arguments here.  One can easily force $Q$ to be parallel 
to all the directions $\dir(\tsdi)$, but we shall only exploit the 
parallelism with $\dir(\tsda)$.

\begin{proof}
Fix $x_0$, $F$, $\tsdi$, and write $A_i = A(x_0,\tsdi)$ for short.  We also
introduce a dummy tube, setting $\tsdnn = \tsda$, $A_n = A_1$.  We choose 
vectors $v_1, \ldots, v_n$ which have the values
$$ v_i = \sd \dir(\tsdi) + O(\delta)$$
for $1 \leq i \leq n$, so that $v_1 = \sd \dir(\tsda)$ and each 
$v_i$ is a distance $\gtrsim \delta$ from the hyperplane spanned by the 
remaining $n-1$ vectors $v_j$. 

The key observation is that for every $\td \in \Td[\tsdi]$,
the set $\td \cap B(x_0,C\sd)$ is essentially constant in
the direction $v_i$.  More precisely, we have the elementary
pointwise estimate 
\be{avv-pt}
 \chi_{B(x_0,C\sd) \cap \td} \lesssim \E_{v_i}( \chi_{B(x_0,2C\sd)
\cap C\td} )
\end{equation}
where $\E_v$ is the averaging operator
$$
\E_v(f(x)) := \int_{|t| \lesssim 1} f(x + v t)\ dt.
$$
Summing this in $\td$ we obtain
\be{ai-av} \chi_{A_i} \lesssim \E_{v_i}(\chi_{\tilde A_i})
\end{equation}

where
$$ \tilde A_i := B(x_0,2C\sd) \cap \bigcup_{\td \in \Td[\tsdi]} C\td.$$

We now claim that $|\tilde A_i| \sim |A_i|$.  Indeed, the lower bound is
trivial, while the upper bound comes from covering $A_i$ by finitely 
overlapping and essentially parallel 
$\delta \times \dots\times\delta\times \sd$ tubes.  

From \eqref{edim-sd} we thus have
\be{tai-size}
|\tilde A_i| \lessapprox \delta^{-C\se} \delta^{n-d} \sd^d.
\end{equation}

To utilize \eqref{ai-av} we invoke

\begin{lemma}\label{three}  Let $v_1, \ldots, v_n$, 
be any $n$ linearly independent vectors in $\R^n$, and
let $F$ be a subset of $\R^n$.
Then for any functions
$f_1, \ldots, f_n$ on $\R^n$, we have
$$ \int_{F}
\prod_{i=1}^n \E_{v_i}(f_i)
\lesssim \sup_P (\frac{|F \cap P|}{|P|})^{1/n}
\prod_{i=1}^n \|f_i\|_{n},$$
where $c > 0$ is an absolute constant, and $P$ ranges over all
parallelepipeds with edge vectors $v_1, \ldots, v_n$.
\end{lemma}

We remark that the $n=3$ version of this lemma was proven in \cite{KLT}.

\begin{proof}  We begin with some reductions.  The statement of the lemma 
is invariant under
affine transformations, so we may rescale $v_i = e_i$, where $e_i$ are the
standard basis of $\R^n$.  It suffices to show that
$$ \int_P \chi_F
\prod_{i=1}^n \E_{e_i}(f_i)
\lesssim (|F \cap P|)^{1/n}
\prod_{i=1}^n \|f_i\|_{L^n(CP)},$$
for all unit cubes $P$, since the claim follows by summing over
a partition of $\R^n$ and using H\"older's inequality.  We may
assume that $P$ is centered at the origin, that $F \subset P$,
and that $f_i$ are supported on $CP$.

By another application of H\"older's inequality, it thus suffices to show that
$$
\| \prod_{i=1}^n \E_{e_i}(f_i) \|_{L^{n/(n-1)}(\R^n)} \lesssim 
\prod_{i=1}^n \|f_i\|_n
$$
for all functions $f_i$ on $CP$.  In fact we shall show the more general statement
\be{n-est}
\| \prod_{i=1}^n \E_{e_i}(f_i) \|_{L^{p/(n-1)}(\R^n)} \lesssim 
\prod_{i=1}^n \|f_i\|_p
\end{equation}
for all $1 \leq p \leq \infty$.

We prove \eqref{n-est} by induction on $n$.  When $n=1$ the claim is clear. 
Now suppose $n>1$, and that \eqref{n-est} has already been proven for 
dimension $n-1$.

For all $x \in \R^n$, write $x = (\underline{x}, \x_n)$, where $\underline{x} \in \R^{n-1}$ and $\x_n$ is the $e_n$ co-ordinate of $x$.

We have the pointwise estimate
$$ \E_{v_n}(f_n)(x) \lesssim \f_1(\underline{x}),$$
where 
$$ \f_n(\underline{x}) = \int_{|x_n| \lesssim 1} f_n(\underline{x},x_n)\ 
dx_n.$$
We can then estimate the left-hand side of \eqref{n-est} by
$$ C (\int \| \f_n \prod_{i=1}^{n-1} \E_{e_i}(f_i^{x_n})\|_{L^{p/(n-1)}
(\R^{n-1})}^{(p/(n-1)} \ dx_n)^{(n-1)/p}$$
where $f_i^{x_n}$ is the function on $\R^{n-1}$ defined by
$$ f_i^{x_n}(\underline{x_n}) := f_i(\underline{x}, x_n).$$
By H\"older in $\R^{n-1}$, we can estimate the previous by
$$ C \| \f_n \|_p (\int \| \prod_{i=1}^{n-1} \E_{e_i}(f_i^{x_n})
\|_{L^{p/(n-2)}(\R^{n-1})}^{(p/(n-1)} \ dx_n)^{(n-1)/p}.$$
By the induction hypothesis, we can estimate this by
$$ C \| \f_n \|_p (\int \prod_{i=1}^{n-1} \| f_i^{x_n} \|_p^{p/(n-1)} \ 
dx_n)^{(n-1)/p}.$$
By another H\"older, we may estimate this by
$$ C \| \f_n \|_p \prod_{i=1}^{n-1} (\int \| f_i^{x_n} \|_p^{p} \ 
dx_n)^{1/p}
= C \| \f_n\|_p \prod_{i=1}^{n-1} \|f_i\|_p.$$
From Young's inequality we have $\|\f_n\|_p \lesssim \|f_n\|_p$, and the 
claim follows.
\end{proof}

Combining this estimate with \eqref{ai-av} and \eqref{tai-size} we obtain
\be{a123e}
|F \cap \bigcap_{i=1}^n A_i| \lessapprox
\delta^{-C\se} \delta^{n-d} \sd^d (\sup_P \frac{|F \cap P|}{|P|})^{1/n}
\end{equation}
where $P$ ranges over all parallelepipeds with edge vectors $v_1, \ldots, v_n$.
To complete the proof of \eqref{grainy-intersect} it thus suffices to show that
\be{f-cover}
\sup_P \frac{|F \cap P|}{|P|} \lesssim \sup_Q \frac{|F \cap Q|}{|Q|}
\end{equation}
where $Q$ ranges over all squares parallel to $\dir(\tsda)$.

Let $\pi$ be the hyperplane generated by $v_1, \ldots, v_{n-1}$, and
let $Q_0$ be a square centered at the origin whose long sides lie on $\pi$. 
We can tile $\R^n$ by translates of $Q_0$, and estimate $P$ by the union of 
all the translates of $Q_0$ which intersect $P$.  This will prove 
\eqref{f-cover} provided that
\be{pqp}
P + Q_0 \subset CP
\end{equation}
for some $C$, where $CP$ is the dilate of $P$ by $C$ around the center of $P$.

To show this, suppose for contradiction that \eqref{pqp} failed.  By 
translation we may assume that $P$ is centered at the origin. 
The failure of \eqref{pqp} then implies that $Q_0$ is not completely
contained inside $CP$. Since $CP$ is convex and symmetric around the
origin, this implies by duality that $CP$ is contained in some slab
$\{x\in \R^n:\ |x\cdot v|\leq 1\}$ for some $v$ outside of
$Q_0^*$, the dual box of $Q_0$.

The dual box $Q_0^*$ has dimensions $\delta^{-1} \times (\dt)^{-1} \times 
\ldots \times (\dt)^{-1}$, is centered at the origin, and has its short 
sides on $\pi$.    Split $v = v_\pi + v_{\pi^\perp}$, where $v_\pi, 
v_{\pi^\perp}$ are the orthogonal projections onto $\pi$ and the orthogonal 
complement of $\pi$ respectively.  Since $v \neq Q_0^*$, we either have 
$|v_{\pi^\perp}| > \delta^{-1}$, or $|v_{\pi^\perp}| \leq \delta^{-1}$
and $|v_\pi| > (\dt)^{-1}$.  In the former case $v_1, \ldots, v_n$ lie
within a $C^{-1}\delta$-neighbourhood of the hyperplane orthogonal to $v$, 
contradicting the choice of the $v_i$.  In the latter case $v_1, \ldots,
v_{n-1}$ lie in the $C^{-1}\dt$-neighbourhood of an $n-2$-dimensional subspace 
of $\pi$, contradicting the non-degeneracy assumption.  This completes the 
proof of \eqref{pqp}, and the lemma follows.
\end{proof}

\section{The planar case}\label{planar-sec}

In this section we shall make heavy use of the fact that $d = \frac{n+2}{2}$,
and so shall perform this substitution throughout the section.  Also, we 
shall be working almost exclusively at scale $\sd$ rather than at $\delta$, 
and so we shall write all of our bounds in terms of $\sd$ rather than $\delta$.

Lemma \ref{triple} gives some control of the intersections of the sets 
$A(x_0,\tsdi)$ provided that the tubes $\tsdi$ are not coplanar.  In order 
to use this Lemma we must ensure that the tubes $\tsd$ which pass through 
a given point $x$ are not too concentrated in a low dimensional space.  
This motivates

\begin{definition}\label{plany-def}  A point $x \in \R^n$ is said to be \emph{degenerate} if there exists an affine subspace $V(x) \subset \R^n$ 
containing $x$ of dimension $n-2$ such that 
\be{clump}
\# \{ \tsd \in \Tsd(x): \tsd \subset \N_{N^C \sd}(V(x)) \} 
\gtrapprox N^{-1} \rho^{-(n-2)/2}.
\end{equation}
If $x$ is not degenerate, we call it \emph{non-degenerate}.
\end{definition}

This bound should be compared to \eqref{many-fat}; in the language of 
\cite{KLT}, it is akin to saying that the set $E_\rho$ is not ``plany'' 
with codimension 2.  The main result of this section is

\begin{proposition}\label{planar}
We have
\be{planar-ok} \allt{\tsd,x} \tsd \in \Tsd, x \in \tsd : 
x \hbox{ is non-degenerate}.
\end{equation}
\end{proposition}

This part of the argument will have a different flavor to the rest of the 
paper.  We remark that the methods used to prove this proposition are not 
used elsewhere in the argument.

Before we begin the rigorous proof of Proposition \ref{planar}, we first 
give an informal argument.  If \eqref{planar-ok} failed, then for most 
points $x \in E_\sd$, a large fraction of the tubes $\tsd$ that pass 
through $x$ will lie near an $n-2$-dimensional space $V(x)$.

Let $\tsda$, $\tsdb$ be a generic pair of tubes intersecting at a point 
$x_0$. Consider the associated ``hairbrush'' associated to $\tsdb$
$$ Brush(\tsdb) := \bigcup_{\tsd \in \Tsd: \tsd \cap \tsdb \neq \emptyset} \tsd.$$
The arguments in Wolff \cite{W1} show that such sets have measure 
\be{brush-bound}
|Brush(\tsdb)| \gtrapprox \sd^{(n-2)/2}.
\end{equation}
Now let $x_1$ be a generic point on $\tsda$, and consider the ``fan'' 
associated to $x_1$
$$ Fan(x_1) := \bigcup_{\tsd \in \Tsd: x_1 \in \tsd} \tsd.$$
The set $Fan(x_1)$ is mostly contained in a small neighbourhood of 
$V(x_1)$.  In particular, $x_0$ should be in this neighborhood.  

Let $\pi$ be the hyperplane spanned by $V(x_1)$ and $\dir(\tsdb)$.  From 
the above considerations we see that $\tsdb$ and $Fan(x_1)$ are both in a 
small neighbourhood of $\pi$.  Thus, we expect that the only tubes in 
$Brush(\tsdb)$ which intersect $Fan(x_1)$ are those which lie in a small 
neighbourhood of $\pi$.  However, an argument from \cite{W2}, \cite{LT1} 
shows that very few tubes in $Brush(\tsdb)$ can be compressed into such a 
small region.  This means that $Brush(\tsdb)$ has a small intersection with
$Fan(x_1)$.  Letting $x_1$ range over all points in $\tsda$, we thus 
conclude that $Brush(\tsda)$ and $Brush(\tsdb)$ have small intersection.  
This can be used together with \eqref{brush-bound} to contradict 
\eqref{emush-sd}.

\begin{figure}[htbp] \centering
\ \psfig{figure=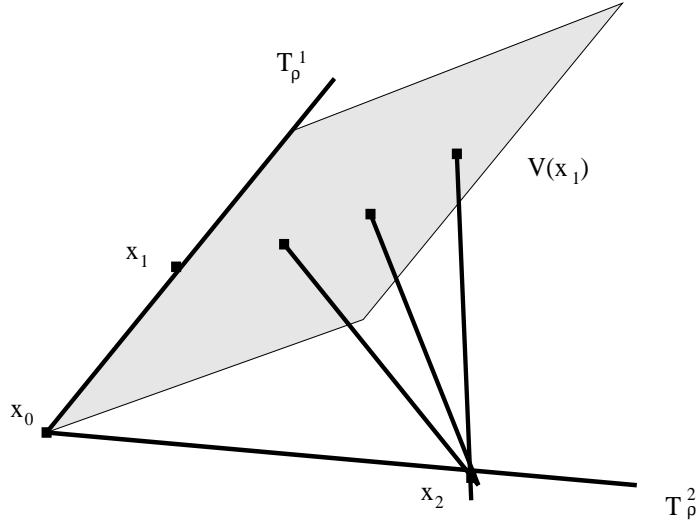,height=2.7in,width=3.6in}
\caption{The only tubes in $Brush (T_\rho^2)$ that intersect $V(x_1)$
are those in a small neighbourhood of the affine subspace spanned
by $V(x_1)$ and dir$(T_\rho^2)$. For clarity, $\rho$-tubes are represented
as lines.
 }
        \label{fig:planar}
        \end{figure}

In order for the above argument to work, one needs a certain amount of 
separation between the various objects under discussion (e.g. one wants 
$|x_1 - x_0|$ and $\angle \dir(\tsdb), V(x_1)$ to be large).  This requires 
a certain amount of technical maneuvering in the rigorous proof, which we 
now begin.

\begin{proof}
Suppose for contradiction that \eqref{planar-ok} failed.  Then we have
$$ | \{ (x,\tsd): \tsd \in \Tsd, x \in \tsd, x \hbox{ degenerate} \}|
\gtrapprox \delta^{c\se} | \{ (x,\tsd): \tsd \in \Tsd, x \in \tsd \}|$$
for some $c > 0$.  From \eqref{p2-ok} we thus have
$$ | \{ (x,\tsd): \tsd \in \Tsd, x \in \tsd, P_2(x), x \hbox{ degenerate} \}|$$
$$\gtrapprox \delta^{c\se} | \{ (x,\tsd): \tsd \in \Tsd, x \in \tsd \}|.$$
From \eqref{tsd-card}, \eqref{filled} the right-hand side is 
$\approx \delta^{c\se}$.  From \eqref{many-fat} we thus have
\be{deg-size}
| \{ (x,\tsd): \tsd \in \Tsd, x \in \tsd, P_2(x), x \hbox{ degenerate} \} 
| \gtrapprox \delta^{C\se} \sd^{(n-2)/2}.
\end{equation}
For every degenerate $x$, let $V(x)$ be an affine subspace satisfying 
\eqref{clump}; one can easily ensure that $V$ is a measurable function.  
Let $\Omega$ denote the set
\be{omega-def}
\Omega := \{ (x, \tsd) \in E_\sd \times \Tsd: P_2(x), x \hbox{ degenerate}, 
x \in \tsd, \tsd \subset \N_{N^C \sd}(V(x)) \}.
\end{equation}
From \eqref{deg-size} and \eqref{clump} we see that $\Omega$ is very large, in fact
\be{om-size}
|\Omega| \gtrapprox N^{-C}.
\end{equation}

On the other hand, we observe that the $x$-projection of $\Omega$ does not 
concentrate in a thin slab.  (This kind of observation also appears in 
\cite{LT1}, and implicitly in \cite{W2}).

\begin{figure}[htbp] \centering
\ \psfig{figure=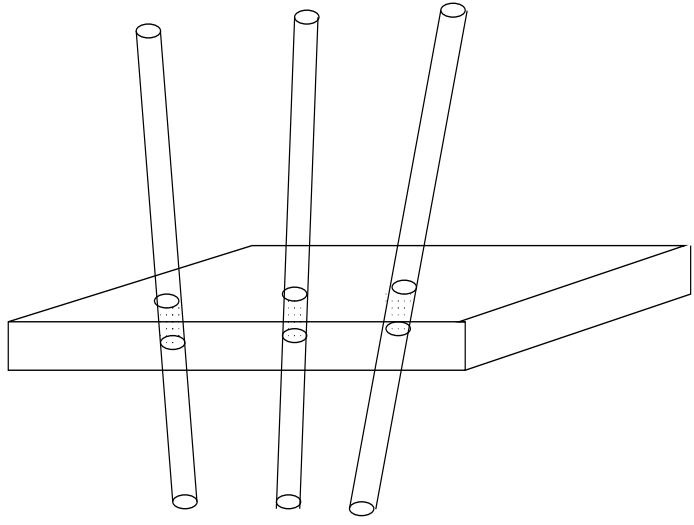,height=2.7in,width=3.6in}
\caption{The set $E_\rho$ (hence the $x$-projection of $\Omega$)
cannot concentrate in thin slabs, since such slabs contain only 
a small fraction of an average tube in $\T_\rho$.
 }
        \label{fig:slab}
        \end{figure}

\begin{lemma}\label{slab}
If $\sd \leq \theta \leq 1$, and $\pi$ is a hyperplane in $\R^n$, then
\be{slab-eq}
| \{ x \in \N_\theta(\pi): (x,\tsd) \in \Omega \hbox{ for some } \tsd \in \Tsd \} |
\lessapprox \delta^{-C\se} \theta \sd^{(n-2)/2}.
\end{equation}
\end{lemma}

Note that the $\theta$ factor on the right-hand side of \eqref{slab-eq} 
gives an improvement over the trivial estimate coming from \eqref{emush-sd}.

\begin{proof}
Let $X$ denote the set on the left-hand side of \eqref{slab-eq}.  From 
\eqref{omega-def} and \eqref{many-fat} we see that
$$ \sum_{\tsd \in \Tsd} \chi_{\tsd}(x) \gtrapprox \delta^{C\se} \sd^{-(n-2)/2}$$
for all $x \in X$.  Integrating this on $X$, we obtain
$$ \sum_{\tsd \in \Tsd} |\tsd \cap X| 
\gtrapprox \delta^{C\se} |X| \sd^{-(n-2)/2}.$$
From elementary geometry we have
$$ |\tsd \cap X| \leq |\tsd \cap \N_\theta(\pi)| \lesssim \frac{\sd^{n-1} \theta}{\theta + \angle(\dir(\tsd),\pi)}.$$
Summing this in $\tsd$, using the $\sd$-separated nature of the directions $\dir(\tsd)$, one obtains
$$ \sum_{\tsd \in \Tsd} |\tsd \cap X| \lessapprox \theta.$$
The claim follows by combining the above estimates.
\end{proof}

The idea is to derive a contradiction by interacting \eqref{om-size} 
with \eqref{slab-eq}.

Let $\T \subset \Tsd$ denote those tubes $\tsd \in \Tsd$ such that
\be{T-def}
| X[\tsd] | \gtrapprox N^{-C} \sd^{n-1},
\end{equation}
where $X[\tsd] \subset \tsd$ is the set
$$
X[\tsd] := \{ x \in \tsd: (x,\tsd) \in \Omega \}.
$$
If the constant $C$ in \eqref{T-def} is chosen sufficiently large, then we 
see from \eqref{tsd-card}, \eqref{T-def}, \eqref{om-size} that
$$ |\{ (x,\tsd) \in \Omega: \tsd \not \in \T \}| \leq \frac{1}{2} |\Omega|$$
and so by \eqref{om-size} again we have
$$ |\{ (x,\tsd) \in \Omega: \tsd \in \T \}| \gtrapprox N^{-C}.$$
In particular, from \eqref{ts-volume} we have
\be{T-card}
\# \T \gtrapprox N^{-C} \sd^{-(n-1)}.
\end{equation}

We now use \eqref{slab-eq} and the ``hairbrush'' argument of Wolff \cite{W1}
to show that the tubes in a hairbrush in $\T$ cannot concentrate in a thin slab.

\begin{lemma}\label{tube-thin}
If $\sd \leq \theta \leq 1$, $\pi$ is a hyperplane in $\R^n$, and 
$\tsdo \in \T$, then
\be{tube-slab-eq}
|\{ (x,\tsd) \in E_\sd \times \T: \tsd \subset \N_\theta(\pi), (x,\tsd), (x, \tsdo) \in \Omega \}|
\lessapprox N^C \theta^c \sd^{n/2}.
\end{equation}
Here $c>0$ is an absolute constant depending only on $n$.
\end{lemma}

As with \eqref{slab-eq}, the key point of \eqref{tube-slab-eq} is that it contains the decay $\theta^c$.

\begin{proof}
We first dispose of the portion where $\angle \tsd,\tsdo \lessapprox 
\theta^{c_0}$, where $c_0 >0$ is some small constant.  In this case we note 
that every $x$ which contributes to \eqref{tube-slab-eq} must satisfy 
$P_2(x)$ and hence \eqref{no-cluster}.  In particular, each $x$ can 
contribute at most $$ \lessapprox \delta^{-C\se} \theta^{c c_0} 
\sd^{-(n-2)/2}$$ tubes $\tsd$ to \eqref{tube-slab-eq}.   Since $x \in 
\tsdo$, the claim then follows from \eqref{ts-volume} and Fubini's theorem.

Now consider the contribution when 
\be{tsd-sep}
\angle \tsd, \tsdo \gtrapprox \theta^{c_0}.
\end{equation}
  Let $\T'$ denote all the tubes in $\T$ which contribute to this portion 
of \eqref{tube-slab-eq}.  By elementary geometry, each $\tsd \in \T'$ 
contributes a set of measure $O(\theta^{-c_0} \sd^n)$ to 
\eqref{tube-slab-eq}.  Thus it suffices to show that
\be{t-targ}
\# \T' \lessapprox N^C \theta^{c+c_0} \sd^{-n/2}.
\end{equation}
For each $\tsd \in \T'$, let $X'[\tsd] \subset X[\tsd]$ denote the set
$$ X'[\tsd] = \{ x \in X[\tsd]: \dist(x,\tsdo) \gtrapprox N^{-C} \theta^{Cc_0} \}.$$
From \eqref{T-def}, \eqref{tsd-sep}, and elementary geometry we see that
$$ |X'[\tsd]| \gtrapprox N^{-C} \sd^{n-1}$$
if the constants are chosen appropriately.  Thus we have
$$ \| \sum_{\tsd \in \T'} \chi_{X'[\tsd]} \|_1 \gtrapprox N^{-C} \sd^{n-1} \# T'.$$
On the other hand, the function $\sum_{\tsd \in \T'} \chi_{X'[\tsd]}$ is 
supported on the set in \eqref{slab-eq}.  From Cauchy-Schwarz we thus have
\be{cord}
\| \sum_{\tsd \in \T'} \chi_{X'[\tsd]} \|_2^2 \gtrapprox N^{-C}\delta^{n-1} (\# T')^2 \theta^{-1} \sd^{-(n-2)/2}.
\end{equation}
We now use a C\'ordoba-style argument.  We can expand the left-hand side as
$$ \sum_{\tsda, \tsdb \in T'} |X'[\tsda] \cap X'[\tsdb]|.$$
We split this sum dyadically based on the angle between $\tsda$ and $\tsdb$:
$$ \sum_{\sd \lesssim 2^{-k} \lesssim 1} \sum_{\tsda, \tsdb \in T': \sd + \angle(\tsda, \tsdb) \sim 2^{-k}} |X'[\tsda] \cap X'[\tsdb]|.$$
Fix $\tsda$.  From elementary geometry, a tube $\tsdb$ can only contribute 
to the sum if it lies within $N^C \theta^{-Cc_0} \sd$ of the 2-dimensional 
plane generated by $\tsdo$ and $\tsda$, and even then the contribution is 
$O(2^k \sd^n)$.  From the $\sd$-separated nature of the tubes $\tsdb$ we 
thus see that there are only $O(2^{-k} \sd^{-1})$ tubes $\tsdb$ which 
contribute to the inner sum.  Combining these observations we thus have
$$
\hbox{LHS of \eqref{cord}} \lessapprox \sum_{\sd \lesssim 2^{-k} \lesssim 1}
2^{-k} \sd^{-1} (\# T') 2^k \sd^n \lessapprox \sd^{n-1} \# T'.
$$
Inserting this into \eqref{cord} and doing some algebra we obtain 
\eqref{t-targ} as desired, if $c_0$ is chosen sufficiently small.
\end{proof}

We now use \eqref{om-size} and the low dimension of the $V(x)$ to 
contradict \eqref{tube-slab-eq}.  

For each $\tsda \in \T$, we have
$$ |\{ (x_0,x_1): x_0, x_1 \in X[\tsda] \}| \gtrapprox N^{-C} \sd^{2(n-1)}$$
by \eqref{T-def}.  We may clearly improve this to
$$ |\{ (x_0,x_1): x_0, x_1 \in X[\tsda]; |x_0 - x_1| \gtrapprox N^{-C} \}| \gtrapprox N^{-C} \sd^{2(n-1)}$$
for appropriate choice of constants.
Summing this over all $\tsda$ and using \eqref{T-card} we obtain
$$ |\{ (x_0,x_1,\tsda): \tsda \in \T; x_0, x_1 \in X[\tsda]; |x_0 - x_1| \gtrapprox N^{-C} \}| \gtrapprox N^{-C} \sd^{n-1}.$$
We rewrite this as
$$ \int_{E_\sd} |\{ (x_1,\tsda): \tsda \in \T; x_0, x_1 \in X[\tsda]; 
|x_0 - x_1| \gtrapprox N^{-C} \}|\ dx_0 \gtrapprox N^{-C} \sd^{n-1}.$$
From \eqref{emush-sd} and Cauchy-Schwarz we thus have
$$ \int_{E_\sd} |\{ (x_1,\tsda): \tsda \in \T; x_0, x_1 \in X[\tsda];
|x_0 - x_1| \gtrapprox N^{-C} \}|^2\ dx_0$$
$$\gtrapprox N^{-C} \sd^{2(n-1)} \sd^{-(n-2)/2}.$$
We write this out as
\be{s1}
|\Sigma|
\gtrapprox N^{-C} \sd^{(3n-2)/2}
\end{equation}
where
$$ \Sigma := \{ (x_0,x_1,x_2,\tsda,\tsdb): \tsda,\tsdb \in \T; x_0,x_1 
\in X[\tsda]; x_0,x_2 \in X[\tsdb];$$
$$|x_0 - x_1|, |x_0 - x_2| \gtrapprox N^{-C} \}.$$
We now claim that

\begin{lemma}\label{s2-lemma}
For any $\sdt \leq \theta \leq 1$, we have
\be{s2}
|\{ (x_0,x_1,x_2,\tsda,\tsdb) \in \Sigma: \angle(\tsdb, V(x_1)) \leq \theta \}| \lessapprox N^C \theta^c \sd^{(3n-2)/2}.
\end{equation}
\end{lemma}

\begin{proof}
From \eqref{emush-sd} and \eqref{tsd-card} it suffices to show that
$$ 
|\{ (x_0,x_2,\tsdb): (x_0,x_1,x_2,\tsda,\tsdb) \in \Sigma; 
\angle(\tsdb, V(x_1)) \leq \theta \}| \lessapprox N^C \theta^c \sd^{(3n-2)/2}$$
for every $x_1$, $\tsda$.  For fixed $x_0$, $\tsdb$, the set of $x_2$ which 
can contribute is $O(\delta^{(n-1)/2})$ by \eqref{ts-volume}.  So it 
suffices to show that
$$ 
|\{ (x_0,\tsdb): (x_0, \tsda), (x_0,\tsdb) \in \Omega; \angle(\tsdb, 
V(x_1)) \leq \theta \}| \lessapprox N^C \theta^c \sd^{n/2}$$
for all non-degenerate $x_1$ and $\tsda \in \Tsd$ such that $\tsda \in \N_{C\sdt}(V(x_1))$.

Fix $x_1$, $\tsda$.  From our assumptions on $\tsdb$ we see that
$\tsdb \in \N_\theta(V(x_1))$.  The claim then follows from 
\eqref{tube-slab-eq} (since $V(x_1)$ can of course be embedded in a hyperplane).
\end{proof}

Combining \eqref{s1}, \eqref{s2} we see that
\be{s3}
|\{ (x_0,x_1,x_2,\tsda,\tsdb) \in \Sigma: \angle(\tsdb, V(x_1)) \gtrapprox N^{-C} \}| \gtrapprox N^{-C} \sd^{(3n-2)/2}
\end{equation}
for appropriate choices of constants.  We rewrite this as
$$ 
\int_{E_\sd}\int_{E_\sd} |\{ (x_0,\tsda,\tsdb): (x_0,x_1,x_2,\tsda,\tsdb)
 \in \Sigma; \angle(\tsdb, V(x_1)) \gtrapprox N^{-C} \}|\ dx_1 dx_2$$
$$\gtrapprox N^{-C} \sd^{(3n-2)/2}.
$$
Using \eqref{emush-sd} and Cauchy-Schwarz as before we thus have
$$ 
\int_{E_\sd}\int_{E_\sd} |\{ (x_0,\tsda,\tsdb): (x_0,x_1,x_2,\tsda,\tsdb)
 \in \Sigma; \angle(\tsdb, V(x_1)) \gtrapprox N^{-C} \}|^2\ dx_1 dx_2$$
$$\gtrapprox N^{-C} \sd^{3n-2} \sd^{-(n-2)},
$$
which we write out as
\be{e1}
\begin{split}
|\{(x_0,x_1,x_2,x_3,&\tsda,\tsdb,\tsdc,\tsdd): 
 (x_0,x_1,x_2,\tsda,\tsdb), (x_3,x_1,x_2,\tsdc,\tsdd) \in \Sigma;\\
&\angle(\tsdb, V(x_1)), \angle(\tsdd,V(x_1)) \gtrapprox N^{-C} \}|
\gtrapprox N^{-C} \sd^{2n}.
\end{split}
\end{equation}

\begin{figure}[htbp] \centering
\ \psfig{figure=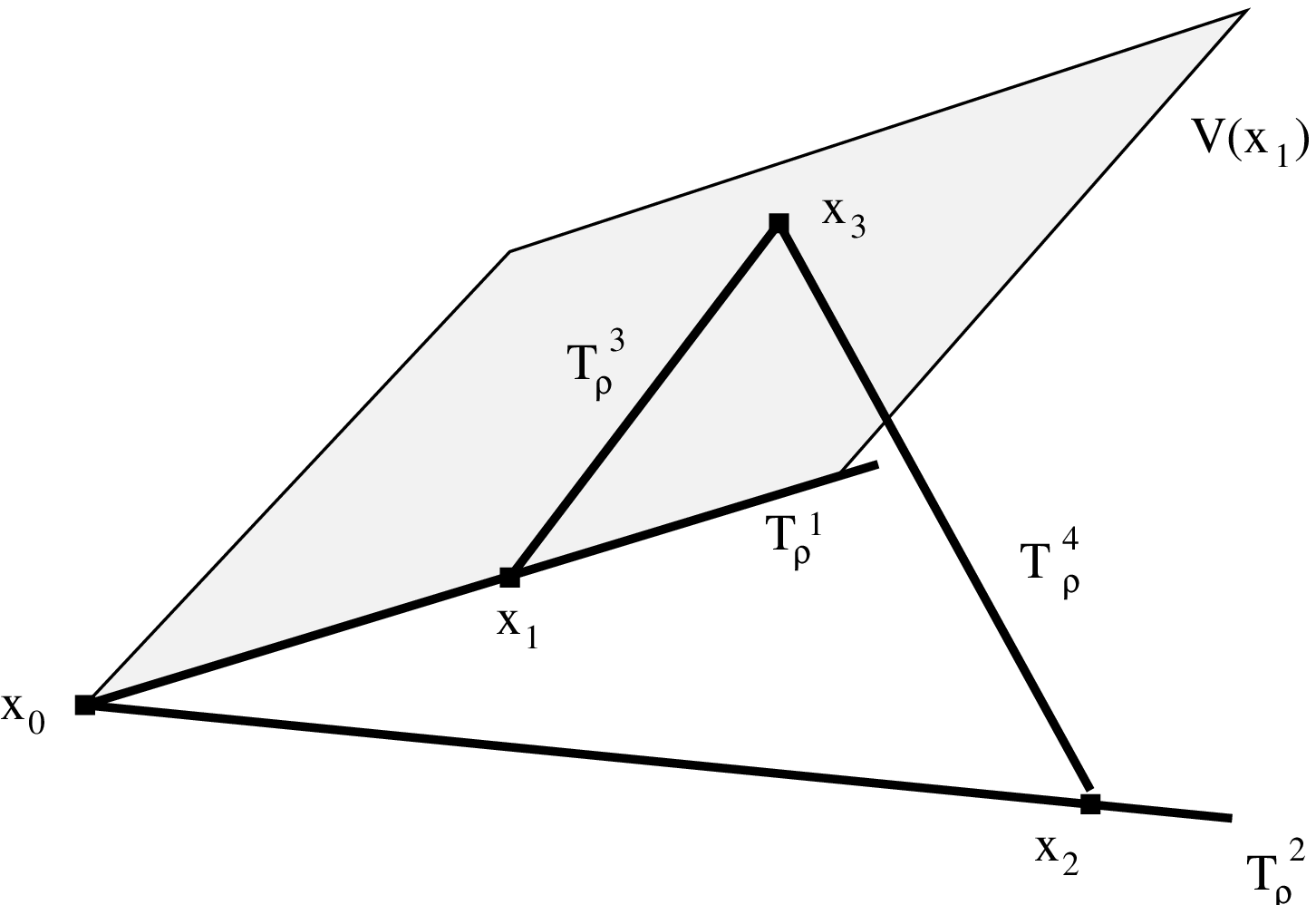,height=2.7in,width=3.6in}
\caption{The $\rho$-tubes $\tsda,\tsdb,\tsdc,\tsdd$ as in \eqref{e1}.
 }
        \label{fig:nonsticky}
        \end{figure}

We now find an upper bound for the left-hand side of \eqref{e1} which will 
achieve the desired contradiction.  The key lemma is

\begin{lemma}\label{e2-lemma}
For each $x_0, x_1, \tsda, \tsdb$, we have
\be{e2}
\begin{split}
|\{(x_2,x_3,\tsdc,\tsdd):&
 (x_0,x_1,x_2,\tsda,\tsdb), (x_3,x_1,x_2,\tsdc,\tsdd) \in \Sigma; \\
&\angle(\tsdb, V(x_1)), \angle(\tsdd,V(x_1)) \gtrapprox N^{-C} \}|
\lessapprox N^C \delta^c \sd^{3n/2}.
\end{split}
\end{equation}
for some absolute constant $c>0$. 
\end{lemma}

\begin{proof}
Fix $x_0, x_1, \tsda, \tsdb$, and let $\pi$ be the hyperplane containing 
$V(x_1)$ and parallel to $\dir(\tsdb)$.  (Note that this hyperplane is well 
defined thanks to the condition $\angle(\tsdb, V(x_1)) \gtrapprox N^{-C}$). 
In order for $(x_2,x_3,\tsdc,\tsdd)$ to contribute to \eqref{e2}, we must have
$$\tsda, \tsdc \subset \N_{N^C \sd}(V(x_1)) \subset \N_{N^C \sd}(\pi).$$
In particular, we have
$$x_0 \in \tsda \subset \N_{N^C \sd}(V(x_1)) \subset \N_{N^C \sd}(\pi).$$
Since $x_0 \in \tsdb$ and $\dir(\tsdb)$ is parallel to $\pi$, we thus have
$$ \tsdb \subset \N_{N^C \sd}(\pi).$$
Since $x_2 \in \tsdb, x_3 \in \tsdc$, we conclude that
$$ x_2, x_3 \in \N_{N^C \sd}(\pi).$$
Since $|x_2 - x_3| \gtrapprox N^{-C}$ and $x_2,x_3 \in \tsdd$, we thus conclude that
$$
\tsdd \subset \N_{N^C \sd}(\pi).
$$
Also, we have
$$
x_3 \in \tsdc \cap \tsdd \subset \N_{C\sdt}(V(x_1)) \cap \tsdd.$$
Using that $\angle(\tsdb,V(x_1)) \gtrapprox N^{-C}$, we see from elementary 
geometry that for fixed $x_2$, $\tsdd$ the set of all possible $x_3$ which 
contribute is contained in a set of measure $\lessapprox N^C \sd^{n}$.  
Also, for fixed $x_2$, $\tsdd$, $x_3$ there is at most $\lessapprox N^C$ 
possible tubes $\tsdc$ which contribute, thanks to the separation condition 
$|x_0 - x_2| \gtrapprox N^{-C}$.  Combining all these observations
we can thus estimate the left-hand side of \eqref{e2} by
$$
\lessapprox N^C \sd^{n} | \{(x_2,\tsdd): \tsdd \subset \N_{N^C \sd}(V(x_1)); (x_2, \tsdb), (x_2,\tsdd) \in \Omega \} |.
$$
The claim then follows from \eqref{tube-slab-eq}.
\end{proof}

In light of \eqref{e2} we may estimate the left-hand side of \eqref{e1} by
$$ \lessapprox N^C \delta^c \sd^{3n/2}
|\{ (x_0,x_1,\tsda,\tsdb): \tsda,\tsdb \in \T; x_0, x_1 \in X[\tsda]; x_0 \in X[\tsdb] \}|.$$
In order for $x_0$ to contribute to the above, $P_2(x_0)$ and thus 
\eqref{many-fat} must hold.  In particular, there are at most $\lessapprox
\delta^{-C\se} \sd^{-(n-2)/2}$ tubes $\tsdb$ which can contribute for each 
$x_0$.  We thus  have
$$ \hbox{LHS of \eqref{e1}} \lessapprox N^C \delta^c \sd^{3n/2}
\sd^{-(n-2)/2} 
|\{ (x_0,x_1,\tsda): \tsda \in \Tsd: x_0, x_1 \in \tsda \}|.$$
From \eqref{tsd-card} and \eqref{ts-volume} we have
$$ |\{ (x_0,x_1,\tsda): \tsda \in \Tsd; x_0, x_1 \in \tsda \}| \lessapprox
\sd^{-(n-1)} \sd^{n-1} \sd^{n-1}.$$
Combining these two estimates together we obtain a contradiction to 
\eqref{e1}, if $\eps$ and then $\delta$ is chosen sufficiently small, and 
the constant $K$ used to define $N$ was chosen sufficiently large so
that $\delta^c\ll N^{-C}$.
\end{proof}

\section{Graininess}\label{grainy-sec}

From Propositions \ref{p3-refine} and \ref{planar} we have
\be{p3-planar-ok} 
\allt{\tsd,x} \tsd \in \Tsd, x \in \tsd : x \hbox{ is non-degenerate}, P_2(x).
\end{equation}
Comparing this with Lemma \ref{triple} we thus expect 
\eqref{grainy-intersect} to happen quite often.  In order to exploit this, 
we shall split $E_\delta$ into a portion which is covered by a small number 
of squares, plus a remainder set $F$ for which we have some control on the 
quantity $(\sup_Q \frac{|F \cap Q|}{|Q|})^{1/n}$.  It turns out that such 
control is essentially automatic for $n > 4$, and for $n=4$ it holds outside
of a small number of squares at each $\dt$-ball.  More precisely, we have

\begin{lemma}\label{grainy-decomp}
Let $B$ be a ball of radius $\dt$ such that
\be{b-bound}
|E_\delta \cap CB| \lessapprox \delta^{-C\se} \delta^{(n-2)/2} (\dt)^{(n+2)/2}
\end{equation}
and let $\omega$ be a direction.  Then we can find a collection $\Q[B, \omega]$ of squares parallel to $\omega$ of cardinality 
\be{q-card}
\# \Q[B, \omega] \lessapprox \delta^{-C\se}
\end{equation}
such that
\be{f-bound}
\sup_Q \frac{|F[B,\omega] \cap Q|}{|Q|} \lessapprox \delta^{C\se}
\end{equation}
where $Q$ ranges over all squares parallel to $\omega$, and
\be{fb-def}
F[B,\omega] := (E_\delta \cap B) \backslash \bigcup_{Q \in \Q[B, \omega]} CQ.
\end{equation}
If $n > 4$ then we can take $\Q[B,\omega]$ to be the empty set.
\end{lemma}

Note that the bound \eqref{b-bound} is consistent with \eqref{edim}.

\begin{proof}
When $n>4$ the claim is trivial with $\Q[B,\omega]$ empty if $\eps$ is 
sufficiently small, since
$$ |F[B,\omega] \cap Q| \leq |E_\delta \cap B|
\lessapprox \delta^{-C\se} \delta^{(n-2)/2} (\dt)^{(n+2)/2}$$
while
$$ |Q| \sim \delta (\dt)^{n-1}.$$

Now suppose that $n=4$.  We say that two squares $Q, Q'$ are 
\emph{separated} if $Q \not \subset 2Q'$ and $Q' \not \subset 2Q$.  We 
define $\Q[B, \omega]$ to be a maximal pairwise-separated set of squares 
$Q$ which satisfy
$$ |E_\delta \cap B \cap Q| \geq \delta^\se |Q|.$$
It is easy to see that \eqref{f-bound} holds.  To show \eqref{q-card}, we 
take advantage of the known bounds for the Radon transform
$$ Rf(t,\omega) = \int \delta(x\cdot \omega - t) f(x)\ dx$$
which takes functions on $\R^4$ to functions on $\R \times S^3$.  (It is 
also possible to obtain \eqref{q-card} by more elementary means).

From the construction of $\Q[B, \omega]$ we see that
$$ | \{ (t,\omega): R\chi_{\N_{C\delta}(E_\delta \cap B)}(t,\omega) 
\gtrsim \delta^\se (\dt)^3 \} | \gtrapprox
\delta N^{-3} \# \Q[B, \omega].$$
On the other hand, one has the restricted weak-type estimate
$$ \| R \chi_E \|_{4,\infty} \lesssim \| \chi_E \|_{4/3}$$
for all sets $E$ (see \cite{oberlin-stein}).  In particular we have
$$ | \{ (t,\omega): R\chi_{\N_{C\delta}(E_\delta)}(t,\omega) \gtrsim 
\delta^\se \} | \lessapprox \delta^{-4\se} (\dt)^{-12} |\N_{C\delta}
(E_\delta \cap B)|^3.$$
On the other hand, since $E_\delta$ is the union of $\delta$-balls we have
$$ |\N_{C\delta}(E_\delta \cap B)| \lesssim |E_\delta \cap CB|
\lessapprox \delta^{-C\se} \delta (\dt)^3$$
by \eqref{b-bound}.  Combining all these estimates we obtain the result.
\end{proof}

Cover $\R^n$ by a finitely overlapping collection $\B$ of $\dt$-balls.  Let $\B'$ denote the subcollection of those balls $B \in \B$ for which \eqref{b-bound} holds.  For each ball $B$ in $\B'$ and each direction $\omega$, we define $\Q[B,\omega]$, $F[B,\omega]$ as in Lemma \ref{grainy-decomp}.  Define the sets $G[\tsd]$, $F[\tsd]$ by
\be{ftsd-def}
F[\tsd] := \bigcup_{B \in \B'} F[B,\dir(\tsd)]; \quad G[\tsd] := 
\bigcup_{B \in \B'} \bigcup_{Q \in \Q[B,\dir(\tsd)]} Q.
\end{equation}

We now combine \eqref{p3-planar-ok}, Lemma \ref{triple}, and Lemma \ref{grainy-decomp} to obtain

\begin{proposition}\label{grainy}
We have
\be{grainy-ok}
\allt{\tsd,\td, x} \tsd \in \Tsd, \td \in \Td[\tsd], x \in \td: x \in G[\tsd].
\end{equation}
\end{proposition}

This immediately yields the desired contradiction when $n > 4$, since the sets $\Q[B,\dir(\tsd)]$ and hence $G[\tsd]$ are always empty.

\begin{proof}
From \eqref{edim} we have
$$ \allt{\td, x} \td \in \Td, x \in \td: |E_\delta \cap B(x,\dt)| \lessapprox \delta^{-\se} \delta^{(n-2)/2} (\dt)^{(n+2)/2}.$$
In particular, we have
$$
\allt{\tsd,\td, x} \tsd \in \Tsd, \td \in \Td[\tsd], x \in \td: x \in \bigcup_{B \in \B'} B.
$$
Using this, \eqref{fb-def}, and \eqref{ftsd-def}, we find that 
\eqref{grainy-ok} will follow if we can show
$$
\allt{\tsd,\td, x} \tsd \in \Tsd, \td \in \Td[\tsd], x \in \td: x \not \in F[\tsd].
$$
By Proposition \ref{delta-sd}, this is equivalent to
$$
\allt{\tsd,\td,x_0,x} \tsd \in \Tsd, x_0 \in \tsd, \td \in \Td[\tsd], 
x \in \td \cap B(x,\sd): x \not \in F[\tsd].
$$
On the other hand, from \eqref{p3-planar-ok} and Proposition \ref{delta-sd} we have
$$
\allt{\tsd,\td,x_0,x} \tsd \in \Tsd, x_0 \in \tsd, \td \in \Td[\tsd], 
x \in \td \cap B(x,\sd): P_2(x_0), x_0 \hbox{ non-degenerate}.
$$
Also, from \eqref{p1-ok}, \eqref{narrow-local}, and Proposition \ref{delta-sd} we have
$$
\allt{\tsd,\td,x_0,x} \tsd \in \Tsd, x_0 \in \tsd, \td \in \Td[\tsd], 
x \in \td \cap B(x,\sd): \hbox{\eqref{narrow-local} holds}.
$$
It thus suffices to show that
\be{bigmosh}
\begin{split}
|\{ (\tsd,\td,x_0,x): &\tsd \in \Tsd, x_0 \in \tsd, \td \in \Td[\tsd], 
x \in \td \cap B(x_0,\sd), P_2(x_0),\\
& x_0 \hbox{ non-degenerate}, \eqref{narrow-local} \hbox{ holds}, x \in F[\tsd] \}| \\
\lessapprox \delta^{c\se}
|\{ (\tsd,\td,&x_0,x): \tsd \in \Tsd, x_0 \in \tsd, 
 \quad \td \in \Td[\tsd], x \in \td \cap B(x_0,\sd)\}|.
\end{split}
\end{equation}
Consider the right-hand side of \eqref{bigmosh}.  For each $\tsd, \td, x$, 
the set of $x_0$ which contribute has volume $\sim \sd^n$.  For each $\tsd, 
\td$, the set of $x$ which contribute has volume $\sim \delta^{n-1}$ by 
\eqref{ts-volume}.  Finally, the total number of pairs $\tsd, \td$ which 
contribute is $\approx \delta^{1-n}$ by \eqref{tsd-card}, \eqref{filled}.  
So the right-hand side is $\approx \delta^{c\se} \sd^n$.  

Now consider the left-hand side of \eqref{bigmosh}.  Using the sets 
$A(x_0,\tsd)$ defined in \eqref{ax0tsd-def}, we can write this as
$$
\int_{P_2(x_0); x_0 \hbox{ non-degenerate}} \sum_{\tsd \in \Tsd(x_0)} 
(\int_{F[\tsd] \cap A(x_0,\tsd): \hbox{\eqref{narrow-local} holds}} 
\mu_\delta[\tsd](x)\ dx)\ dx_0.$$
By \eqref{narrow-local}, we can estimate this by
$$
\delta^{-\se} \sd^{-(n-d)} \int_{P_2(x_0); x_0 \hbox{ non-degenerate}} 
\sum_{\tsd \in \Tsd(x_0)} |F[\tsd] \cap A(x_0,\tsd)| \ dx_0.$$
We rewrite this as
$$
\delta^{-\se} \sd^{-(n-d)} \int_{P_2(x_0); x_0 \hbox{ non-degenerate}} 
\| \sum_{\tsd \in \Tsd(x_0)} \chi_{F[\tsd] \cap A(x_0,\tsd)} \|_1 \ dx_0.$$
The expression inside the norm is supported inside $B(x_0,\sd) \cap 
E_\delta$, which has measure $\lessapprox \delta^{-C\se} \delta^{n-d} 
\sd^d$ by \eqref{edim-sd}.  By H\"older's inequality, we may therefore 
estimate the above as
\be{mish-mash}
\begin{split}
&\delta^{-C\se} \sd^{-(n-d)} (\delta^{n-d} \sd^d)^{1/(n-1)'} \\
&\times \int_{P_2(x_0); x_0 \hbox{ non-degenerate}} 
\| \sum_{\tsd \in \Tsd(x_0)} \chi_{F[\tsd] \cap A(x_0,\tsd)} \|_{n-1} \ dx_0.
\end{split}
\end{equation}

We now estimate this norm as

\begin{lemma}\label{n-1}
If $P_2(x_0)$ holds and $x_0$ is non-degenerate, then we have
\be{n-1-bound}
\| \sum_{\tsd \in \Tsd(x_0)} \chi_{F[\tsd] \cap A(x_0,\tsd)} \|_{n-1}
\lessapprox \delta^{c\se} \sd^{-(n-d)} (\delta^{n-d} \sd^d)^{1/(n-1)}
\end{equation}
for some absolute constant $c > 0$.
\end{lemma}

\begin{proof}
Fix $x_0$.
Raising both sides of \eqref{n-1-bound} to the $n-1$th power and expanding, 
it suffices to show that
$$ \sum_{\tsda, \ldots, \tsdn \in \Tsd(x_0)} |\bigcap_{i=1}^{n-1}
\chi_{F[\tsdi] \cap A(x_0,\tsdi)}| \lesssim N^{-c} \sd^{-(n-d)(n-1)}
\delta^{n-d} \sd^d.$$
We first deal with the contribution when the tubes $\tsda, \ldots, \tsdn$ 
are not coplanar.  In this case we use Lemma \ref{triple} to estimate the 
above by
$$
\sum_{\tsda, \ldots, \tsdn \in \Tsd(x_0)} \delta^{-C\se} \delta^{n-d} 
\rho^d (\sup_Q \frac{|F[\tsda] \cap Q|}{|Q|})^{1/n},$$
where the supremum is over all squares $Q$ parallel to $\dir(\tsda)$.
By \eqref{f-bound}, \eqref{ftsd-def} and the finite overlap of the balls 
$B$ we have
$$ \sup_Q \frac{|F[\tsda] \cap Q|}{|Q|} \lessapprox \delta^{\se}.$$
The claim then follows from \eqref{many-fat}.

It remains to control the contribution when the tubes $\tsda, \ldots, 
\tsdn$ are coplanar.  In this case we use \eqref{edim-sd} to make the crude 
estimate
$$ |\bigcap_{i=1}^{n-1} \chi_{F[\tsdi] \cap A(x_0,\tsdi)}|
\leq |E_\delta \cap B(x_0,\sd)| \lessapprox \delta^{-C\se} \delta^{n-d} 
\sd^d.$$
By \eqref{many-fat}, it thus suffices to show that
\be{wookie}
\# \{ (\tsda, \ldots, \tsdn) \in \Tsd(x_0)^{n-1}: \tsda, \ldots, \tsdn \hbox{ coplanar} \} \lessapprox N^{-c} (\# \Tsd(x_0))^{n-1}.
\end{equation}
For any $1 \leq k < n-1$ and any tubes $\tsda, \ldots, \tsdk$ we choose an 
$n-2$ dimensional space $V(\tsda, \ldots, \tsdk)$ through $x_0$ which is parallel to 
$\dir(\tsda), \ldots, \dir(\tsdk)$.  This choice of space may not always be 
unique, but we select it in such a way that $V$ is measurable.  From 
elementary geometry we see that if $\tsda, \ldots, \tsdn$ are coplanar, 
then we must have
$$ \tsdka \subset \N_{N^C \sd}(V(\tsda, \ldots, \tsdk))$$
for some $1 \leq k < n-1$.  The claim \eqref{wookie} then follows from the 
non-degeneracy of $x_0$.
\end{proof}

By this lemma, we can estimate \eqref{mish-mash} by
$$
N^{-c}
\sd^{-2(n-d)} (\delta^{n-d} \sd^d) 
\int_{P_2(x_0)} \ dx_0.
$$
Since the integral is clearly bounded by $|E_\sd|$, we can estimate this by $\delta^{c\se} \sd^n$ as desired by \eqref{emush-sd}, if $\eps$ is sufficiently small.
\end{proof}

\section{The grainy four-dimensional case}\label{4D-sec}

We have already proven Theorem \ref{main-thm} when $n > 4$.  Accordingly, 
we shall assume for the remainder of the argument that $n=4$.

The key geometrical observation shall be a ``four-square lemma'', Lemma 
\ref{four-grain}, which places a non-trivial limit on the possible 
directions of $\delta$-tubes which simultaneously pass through four 
separated squares.  (This can be thought of as the four-dimensional 
analogue of the ``three-line lemma'' used in \cite{schlag:kakeya}).

From Lemma \ref{allt-mix} and \eqref{ts-volume} we can rewrite 
\eqref{grainy-ok} as
$$
\allt{\tsd,\td} \tsd \in \Tsd, \td \in \Td[\tsd]: 
(\allt{x} x \in \td: x \in G[\tsd]) 
$$
Clearly, if $\tsd$, $\td$ satisfy
$$
\allt{x} x \in \td: x \in G[\tsd]$$
then they also satisfy
$$
\allt{x_1,x_2,x_3,x_4} x_1,x_2,x_3,x_4 \in \td: 
x_1,x_2,x_3,x_4 \in G[\tsd]$$
(this can either be proved directly, or by iterating Lemma \ref{allt-mix} 
and \eqref{ts-volume}).  Thus we have
$$
\allt{\tsd,\td} \tsd \in \Tsd, \td \in \Td[\tsd]: (\allt{x_1,x_2,x_3,x_4} 
x_1,x_2,x_3,x_4 \in \td: x_1,x_2,x_3,x_4 \in G[\tsd]). 
$$
From Lemma \ref{allt-mix} and \eqref{ts-volume} we can rewrite this as
$$
\allt{\tsd,\td,x_1,x_2,x_3,x_4} \tsd \in \Tsd, \td \in \Td[\tsd], 
x_1,x_2,x_3,x_4 \in \td: x_1,x_2,x_3,x_4 \in G[\tsd]. 
$$
In particular, from \eqref{rarity} we have
$$
| \{ (\tsd,\td,x_1,x_2,x_3,x_4): \ \tsd \in \Tsd, \td \in \Td[\tsd], 
x_1,x_2,x_3,x_4 \in \td \cap G[\tsd] \}|$$ 
$$\sim 
| \{ (\tsd,\td,x_1,x_2,x_3,x_4):
\tsd \in \Tsd, \td \in \Td[\tsd], x_1,x_2,x_3,x_4 \in \td \}|.$$
From \eqref{ts-volume}, \eqref{filled}, \eqref{tsd-card} the right-hand side is
$$ \approx (1/\delta)^{(n-1)}\delta^{4(n-1)}.$$
We therefore have
\be{4a}
\sum_{\tsd \in \Tsd} \sum_{\td \in \Td[\tsd]}
| \{ (x_1,x_2,x_3,x_4): x_1,x_2,x_3,x_4 \in \td \cap G[\tsd] \}|
\approx (1/\delta)^{(n-1)}\delta^{4(n-1)}.
\end{equation}
Let $0 < \theta \ll 1$ be a quantity to be chosen shortly.  From elementary 
geometry we have
$$ | \{ (x_1,x_2,x_3,x_4): x_1,x_2,x_3,x_4 \in \td; |x_1 - x_2| \leq 
\theta \}| \lessapprox \theta \delta^{4(n-1)}$$
and so by \eqref{ts-volume}, \eqref{filled} as before we have
$$
\sum_{\tsd \in \Tsd} \sum_{\td \in \Td[\tsd]}
| \{ (x_1,x_2,x_3,x_4): x_1,x_2,x_3,x_4 \in \td \cap G[\tsd]; |x_1 - x_2| 
\leq \theta \}|$$
$$\lessapprox \theta (1/\delta)^{n-1} \delta^{4(n-1)}.$$
Similarly for permutations of the indices $1,2,3,4$.  Combining these 
estimates with \eqref{4a}, we obtain
\be{chock}
\begin{split}
\sum_{\tsd \in \Tsd} \sum_{\td \in \Td[\tsd]}
| \{ (x_1,x_2,x_3,x_4) &\in (\td \cap G[\tsd])^4: |x_i - x_j| \approx 1 
\hbox{ for } 1 \leq i < j \leq 4\}|\\
&\approx (1/\delta)^{n-1} \delta^{4(n-1)}.
\end{split}
\end{equation}

We now pause to interpose a family of $\dt$-tubes between the 
$\delta$-tubes in $\Td$ and the $\sd$-tubes in $\Tsd$.

\begin{lemma}\label{sticky-dt} 
There exists a family $\Tdt$ of $\dt$-tubes such that
\be{tdt-card}
\# \Tdt \lessapprox \delta^{-\se} (\frac{1}{\dt})^{n-1}
\end{equation}
and such that 
\be{stick}
\# \{ \td \in \Td: \td \not \subset \tdt \hbox{ for all } \tdt \in \Tdt \} 
\lessapprox \delta^{c\se} (\frac{1}{\delta})^{n-1}.
\end{equation}
\end{lemma}

\begin{proof}
Let $\Ecal$ be a maximal $\dt$-separated set of directions, and for each 
$\omega \in E$ let $\Tdt[\omega]$ be a finitely overlapping cover of $\R^n$ 
by $\dt$-tubes with direction $\omega$.  We can arrange matters so that 
every $\td \in \Td$ obeys $\td \subset \tdt$ for some $\omega \in \Ecal$ 
and $\tdt \in \Tdt[\omega]$.

Call a direction $\omega \in \Ecal$ \emph{sticky} if
$$ \# \{ \tdt \in \Tdt[\omega]: \tdt \supset \td \hbox{ for some } \td \in 
\Td \} \leq \delta^{-\se}$$
and define
$$ \Tdt := \{ \tdt: \tdt \in \Tdt[\omega] \hbox{ for some sticky } \omega \in \Ecal;
\tdt \supset \td \hbox{ for some } \td \in \Td \}.$$
Clearly \eqref{tdt-card} holds.  To prove \eqref{stick} it suffices to show that
$$
\# \{ \td \in \Td: \td \subset \tdt \hbox{ for some } \tdt \in \Tdt[\omega] 
\hbox{ and some non-sticky } \omega \in \Ecal \} $$
$$\lessapprox \delta^{c\se} (\frac{1}{\delta})^{n-1}.
$$
Since $\Td$ is direction-separated, each non-sticky direction $\omega$ can 
contribute at most $N^{n-1}$ elements to the above set.  Hence it suffices 
to show that
\be{e-bound}
\# \Ecal' \lessapprox \delta^{c\se} (\frac{1}{\dt})^{n-1}
\end{equation}
where $\Ecal'$ is the set of non-sticky directions.

By construction, for each $\omega \in \Ecal'$ we can find a subset 
$\Tdt'[\omega] \subset \Tdt[\omega]$ of cardinality $\# \Tdt'[\omega] 
\approx \delta^{-\se}$ such that each $\tdt \in \Tdt'[\omega]$ contains at 
least one tube $\td \in \Td$.  Let $\Tdt'$ be the union of all these 
$\Tdt'[\omega]$ as $\omega$ ranges over $\Ecal'$.  By construction, the 
$\Tdt'$ have directional multiplicity $\lessapprox \delta^{-\se}$, and we have
$$ \# \Tdt' \approx \delta^{-\se} \# \Ecal'$$
and
$$ \bigcup_{\tdt \in \Tdt'} \tdt \subset \N_{C\dt}(\bigcup_{\td \in \Td} \td).$$
In particular, from \eqref{emush-sigma} we have
$$ |\bigcup_{\tdt \in \Tdt'} \tdt| \lessapprox (\dt)^{n-d}.$$
On the other hand, from \eqref{max-xray} we have
$$
\| \sum_{\tdt \in \Tdt'} \chi_{\tdt} \|_{d'}
\lessapprox (\dt)^{\frac{d-n}{d}} ((\dt)^{n-1} \# \Ecal')^{\frac{n-2}{n-1}
 + \frac{1}{d(n-1)}}.
$$
From H\"older's inequality we have
$$
\| \sum_{\tdt \in \Tdt'} \chi_{\tdt} \|_{1}
\lessapprox ((\dt)^{n-1} \# \Ecal')^{\frac{n-2}{n-1} + \frac{1}{d(n-1)}}.
$$
However, by \eqref{ts-volume} we have
$$
\| \sum_{\tdt \in \Tdt'} \chi_{\tdt} \|_{1}
\approx \dt^{n-1} \# \Tdt' \approx \delta^{-\se} \dt^{n-1} \# \Ecal'.$$
Combining these two inequalities we obtain \eqref{e-bound} as desired.
\end{proof}

Let $\Tdt$ be as in the above lemma.  Returning to \eqref{chock}, we note 
that each tube $\td \in \Td$ can contribute at most
$$ |\td|^4 \sim \delta^{4(n-1)}$$
to \eqref{chock}.  From this and \eqref{stick} we thus have
\bas
\sum_{\tsd \in \Tsd} &\ \sum_{\td \in \Td[\tsd]: \td \subset \tdt \hbox{ for some } \tdt \in \Tdt}
| \{ (x_1,x_2,x_3,x_4) \in (\td \cap G[\tsd])^4:\\
&|x_i - x_j| \approx 1 \hbox{ for } 1 \leq i < j \leq 4\}| \approx (1/\delta)^{n-1} \delta^{4(n-1)}.
\end{align*}
From \eqref{tdt-card}, there must therefore exist a tube $\tdt \in \Tdt$ such that
\bas
\sum_{\tsd \in \Tsd} \sum_{\td \in \Td[\tsd]: \td \subset \tdt} &
| \{ (x_1,x_2,x_3,x_4) \in (\td \cap G[\tsd])^4: |x_i - x_j| \approx 1 \\
&\hbox{ for } 1 \leq i < j \leq 4\}| \gtrapprox (\sd/\delta)^{n-1}\delta^{4(n-1)}.
\end{align*}
Fix this $\tdt$.  Since $\td$ must be contained in both $\tsd$ and $\tdt$, 
we see from elementary geometry that 
\be{tsd-td}
\dir(\tsd) = \dir(\tdt) + O(\sd).
\end{equation}
  Since the collection $\Tsd$ is direction-separated, we may therefore find 
a tube $\tsd \in \Tsd$ obeying \eqref{tsd-td} such that
\bas
\sum_{\td \in \Td[\tsd]: \td \subset \tdt}
| \{ (x_1,x_2,x_3,x_4) & \in (\td \cap G[\tsd])^4:
|x_i - x_j| \approx 1 \\ 
&\hbox{ for } 1 \leq i < j \leq 4\}| \gtrapprox (\sd/\delta)^{n-1}
\delta^{4(n-1)}.
\end{align*}
Fix this $\tsd$.  Let $\B''$ denote all the balls in $\B'$ which intersect 
$\tdt$.  Note that $B \subset C\tdt$ for all $B \in \B''$.
From \eqref{ftsd-def} we thus have
$$
\sum_{B_1,B_2,B_3,B_4 \in \B'': \dist(B_i,B_j) \approx 1 \hbox{ for } 1 \leq i <  j \leq 4}
\sum_{\td \in \Td[\tsd]: \td \subset \tdt} $$
$$\prod_{i=1}^4 |\td \cap \bigcup_{Q_i \in \Q[B_i,\dir(\tsd)]} Q_i|
\gtrapprox N^{n-1} \delta^{4(n-1)}.
$$
From elementary geometry we have
$$ \# \B'' \lessapprox (\dt)^{-1}.$$
We may therefore find balls $B_1,B_2,B_3,B_4 \in \B''$ satisfying
$$ \dist(B_i,B_j) \approx 1 \hbox{ for } 1 \leq i < j \leq 4 $$
and such that
$$
\sum_{\td \in \Td[\tsd]: \td \subset \tdt} 
\prod_{i=1}^4 |\td \cap \bigcup_{Q_i \in \Q[B_i,\dir(\tsd)]} Q_i|
\gtrapprox N^{n-1} (N\delta^n)^4.
$$
Fix these $B_1$, $B_2$, $B_3$, $B_4$.  From \eqref{q-card} we may thus find squares $Q_i \in \Q[B_i,\dir(\tsd)]$ for $i=1,2,3,4$ such that
$$
\sum_{\td \in \Td[\tsd]: \td \subset \tdt} 
\prod_{i=1}^4 |\td \cap Q_i|
\gtrapprox N^{n-1} (N\delta^n)^4.
$$
Fix $Q_1$, $Q_2$, $Q_3$, $Q_4$; note that $Q_i \in C \tdt$ for $i=1,2,3,4$.  From elementary geometry we have
$$ |\td \cap Q_i| \lessapprox N \delta^n.$$
From the preceding we must therefore have
$$ \# \{ \td \in \Td[\tsd]: \td \subset \tdt; \td \cap Q_i \neq \emptyset 
\hbox{ for } i=1,2,3,4 \} \gtrapprox N^{n-1} = N^3.$$
On the other hand, from the direction-separated nature of the $\td$ we have 
the trivial estimate
\be{no-contra}
\# \{ \td \in \Td[\tsd]: \td \subset \tdt \} \lessapprox N^3.
\end{equation}
These two statements are not quite in contradiction.  However, we can 
obtain the following improvement to \eqref{no-contra}, and this will yield
the desired contradiction.

\begin{lemma}\label{four-grain}
Let $n=4$, $\tsd \in \Tsd$, $\tdt \in \Tdt$ be tubes obeying \eqref{tsd-td},
and let $Q_1, Q_2, Q_3, Q_4$ be four squares in $C\tdt$ parallel to 
$\dir(\tsd)$ such that $\dist(Q_i,Q_j) \approx 1$ for all $1 \leq i < j 
\leq 4$.  Let $\T$ be a collection of direction-separated $\delta$-tubes 
in $\tdt$ such that $T \cap Q_i \neq \emptyset$ for all $T \in \T$, 
$i=1,2,3,4$.  Then
$$ \# \T \lessapprox N^{3 - \frac{1}{4}}.$$
\end{lemma}

The $1/4$ gain is not best possible, but that is irrelevant for our purposes.

To complete the proof of Theorem \ref{main-thm} it only remains to prove 
Lemma \ref{four-grain}.  This we shall do in the next section.

\section{Linear algebra}

We now prove Lemma \ref{four-grain}.  Roughly speaking, this lemma is 
stating that requiring a line to intersect four distinct horizontal 
2-planes must constrain the line to a 2-dimensional set of directions, as 
opposed to the full 3-dimensional set of directions.

By \eqref{tsd-td} we can perturb the $Q_i$ to be parallel to $\dir(\tdt)$ 
rather than $\dir(\tsd)$.  The reader may verify that this has essentially 
no effect on the statement and conclusions of the lemma.  The tube $\tsd$ 
now plays no role and will be ignored.

By an affine transformation we may assume that $\tdt$ is the vertical tube
$$ \tdt = \{ (\underline{x}, \x_n): 0 \leq \x_n \leq 1, |\underline{x}| \leq \dt \}.$$
We may replace each square $Q_i$ by its central horizontal slice
$$ \{ x \in Q_i: \x_n = t_i \}$$
where $t_i$ is the $n$-th co-ordinate of the center of $Q_i$.

If we now apply the non-isotropic rescaling $(\underline{x}, \x_n) \to (\underline{x}/(\dt), \x_n)$ to map $\tdt$ to the unit cylinder, the 
problem now reduces to proving

\begin{lemma}\label{four-grain-alt}
Let $n=4$, and let $t_1, t_2, t_3, t_4$ be four numbers in $[0,1]$ such 
that $|t_i - t_j| \approx 1$ for all $1 \leq i < j \leq 4$.  Let $A_1$, 
$A_2$, $A_3$, $A_4$ be four boxes of dimensions $\frac{C}{N} \times C 
\times C$ in $\R^3$.  Let $\T$ be a collection of direction-separated 
$1/N$-tubes in a bounded region of $\R^4$ such that
$$ T \cap \{t_i \} \times A_i \neq \emptyset$$
for all $T \in \T$, $i = 1,2,3,4$.  Then $\# T \lessapprox N^{3 - \frac{1}{4}}.$
\end{lemma}

\begin{proof}
We can find unit directions $\omega_i \in S^2 \subset \R^3$ and numbers 
$b_i \in \R$ for all $i=1,2,3,4$ such that $b_i = O(1)$ and 
$$ A_i \subset \{x \in \R^3: \omega_i \cdot x = b_i + O(1/N) \},$$
where the dot product is taken in $\R^3$.

Fix $\omega_i$ and $b_i$. Let $T$ be a tube in $\T$.  We can find $x, v 
\in \R^3$ with $|x|, |v| \lesssim 1$ such that 
$$ T \subset \{ (x+vt+O(1/N),t): 0 \leq t \leq 1 \},$$
so in particular we have
\be{xv}
(x+vt_i) \cdot \omega_i = b_i + O(1/N)
\end{equation}
for $i=1,2,3,4$.  Since $\T$ is direction-separated, it thus suffices to
show that the set of all possible velocities $v$ which obey \eqref{xv}
for some $x$ can only support $\lessapprox N^{3-\frac{1}{4}}$
$1/N$-separated values at best.
By linearity, we may assume that $b_i=0$.

We define the \emph{rank} $k$ to be the least integer $k$ such that there
 exist distinct $i_1, \ldots, i_k$ in $\{1,2,3,4\}$ and co-efficients 
$a_1, \ldots, a_k \in \R$ such that
\be{c1}
\max(|a_1|, \ldots, |a_k|) \geq 1
\end{equation}
and
\be{ringo}
|a_1 \omega_{i_1} + \ldots + a_k \omega_{i_k}| \leq N^{-k/4},
\end{equation}
where $0 < c_1 \ll 1$ is an absolute constant to be chosen later.  Since 
the $\omega_i$ live in $\R^3$ and have magnitude 1, we see that the rank 
is well-defined and is either 2, 3, or 4.

Fix $k$ to be the rank, and let $a_1, \ldots, a_k$ be as above.  Clearly
we may normalize so that 
$$ \max(|a_1|, \ldots, |a_k|) = |a_1| = 1.$$
If we multiply \eqref{xv} for $i = i_j$ by $a_j$ for $j = 1, \ldots, k$
and add, we obtain
$$ (x+t_k v) \cdot (a_1 \omega_{i_1} + \ldots + a_k \omega_{i_k})
+ v \cdot (a_1 (t_1-t_k) \omega_{i_1} + \ldots + a_{k-1} (t_{k-1} - t_k) 
\omega_{i_{k-1}}) = O(1/N).$$
From \eqref{ringo} we have
$$ |(x+t_k v) \cdot (a_1 \omega_{i_1} + \ldots + a_k \omega_{i_k})| 
\lessapprox N^{-k/4}$$
whereas by the definition of rank and the fact that $|a_1 (t_1 - t_k)| 
\approx 1$ we have
$$ | a_1 (t_1-t_k) \omega_{i_1} + \ldots + a_{k-1} (t_{k-1} - t_k) 
\omega_{i_{k-1}}| \gtrapprox N^{-(k-1)/4}.$$
Since $O(1/N) = O(N^{-k/4})$, we thus see that $v$ is constrained to lie
in the $\lessapprox N^{-1/4}$-neighbourhood of a hyperplane.  This means
that any $1/N$-separated set of such $v$ can have cardinality at most 
$\lessapprox N^{3 - 1/4}$, as desired.
\end{proof}

\end{document}